\documentclass[11pt,epsfig]{article}
\usepackage{mathrsfs}
\usepackage{amsmath}
\usepackage{amsthm}
\usepackage{epsfig}
\usepackage{graphicx}
\usepackage{hyperref}
\usepackage{amsfonts}
\usepackage{color}
\usepackage[margin=1.2in]{geometry}
\usepackage{setspace}

\newtheoremstyle{thmm}{1.5ex plus 1ex minus .2ex}{1.5ex plus
1ex minus
.2ex}{\rmfamily}{}{\bfseries}{}{1em}{} \theoremstyle{thmm}
\newtheorem{theorem}{Theorem}[section]
\newtheorem{lemma}{Lemma}[section]

\renewenvironment{proof}[1][Proof]{\noindent\textit{#1. }
}{\hfill$\square$}

\allowdisplaybreaks
\renewcommand{\theequation}{\thesection.\arabic{equation}}

\newcommand{\nn}{\nonumber}
\def \endproof{\vrule height8pt width 5pt depth 0pt}
\def\refe#1{(\ref{#1})}

\newcommand{\vertiii}[1]
{{\left\vert\kern-0.25ex\left
\vert\kern-0.25ex\left\vert #1
    \right\vert\kern-0.25ex\right
\vert\kern-0.25ex\right\vert}}

\def\d{\delta}

\def\R{\mathbb{R}}

\def\d{{\rm d}}

\begin{document}

\title{\bf Maximum-norm stability and maximal $\bf L^p$ regularity of FEMs for parabolic equations with Lipschitz continuous coefficients}
\author{
Buyang Li \footnote{
Department
of Mathematics, Nanjing University, Nanjing 210093, Jiangsu, P.R. China.
The author was supported in part by a grant from NSFC
(Grant No. 11301262) \newline
\indent~ Email address: buyangli@nju.edu.cn 
}
}
\date{}
\maketitle

\begin{abstract}
In this paper, we study the semi-discrete Galerkin finite
element method for parabolic equations with
Lipschitz continuous coefficients. We prove
the maximum-norm stability of the semigroup
generated by the corresponding elliptic finite element operator,
and prove the space-time stability of the parabolic
projection onto the finite element space in
$L^\infty( Q_T)$ and
$L^p((0,T);L^q(\Omega))$, $1<p,q<\infty$.
The maximal $L^p$ regularity of
the parabolic finite element equation is also established.

\end{abstract}

\section{Introduction}
\setcounter{equation}{0}  Consider the linear parabolic equation
\begin{align}\label{PDE0}
\left\{
\begin{array}{ll}
\displaystyle
\partial_tu-\sum_{i,j=1}^d\partial_{i}\big(a_{ij}(x)\partial_j
u\big)+c(x)u =f-\sum_{j=1}^d\partial_jg_j
&\mbox{in}~~\Omega\times(0,T),\\[5pt]
\displaystyle \sum_{i,j=1}^da_{ij}(x)n_i\partial_j u=\sum_{j=1}^dg_jn_j
&\mbox{on}~~\partial\Omega\times(0,T),\\[8pt]
u(\cdot,0)=u^0 &\mbox{in}~~\Omega,
\end{array}
\right.
\end{align}
where $\Omega$ is a bounded smooth domain in $\R^d$ $(d\geq 2)$, $T$
is a given positive number, $f$ and ${\bf g}=(g_1,\cdots,g_d)$ are
given functions. The Galerkin finite element method (FEM) for the
above equation seeks $\{u_h(t)\in S_h\}_{t> 0}$
satisfying the parabolic finite element equation:
\begin{align}\label{FEEq0}
\left\{
\begin{array}{ll}
\displaystyle \big(\partial_tu_h,v_h\big)
+\sum_{i,j=1}^d\big(a_{ij}\partial_j
u_h,\partial_iv_h\big)+\big(cu_h,v_h\big)
=(f,v_h)+\sum_{j=1}^d(g_j,\partial_jv_h), ~~\forall~v_h\in
S_h,\\[8pt]
u_h(0)=u^0_h ,
\end{array}
\right.
\end{align}
where $S_h$, $0<h<h_0$, denotes a finite element subspace of
$H^1(\Omega)$ consisting of continuous piecewise polynomials of
degree $r\geq 1$ on certain quasi-uniform triangulations of 
$\Omega$ which fit the boundary exactly.
The coefficients $a_{ij}=a_{ji}$, $i,j=1,\cdots,d$, and $c$ in the
above equations are assumed to satisfy
\begin{align}\label{coeffcond}
\Lambda_1|\xi|^2\leq
\sum_{i,j=1}^da_{ij}(x)\xi_i\xi_j
\leq \Lambda_2|\xi|^2
\quad\mbox{and}\quad
c(x)\geq c_0,\quad\mbox{for}~~x\in\Omega,
\end{align}
for some positive constants $\Lambda_1,\Lambda_2$ and $c_0$.

Define the
elliptic operator $A: H^1(\Omega)\rightarrow H^1(\Omega)'$ and its
finite element approximation $A_h: S_h\rightarrow S_h$
by
\begin{align}
& (Aw,v):=\sum_{i,j=1}^d(a_{ij}\partial_j w,\partial_iv)+(cw,v),
~~~~~\qquad\qquad \forall~w,v\in
H^1(\Omega),\\
 &(A_hw_h,v_h):=\sum_{i,j=1}^d(a_{ij}\partial_j
w_h,\partial_{i}v_h)+(cw_h,v_h),\qquad~ \forall~w_h,v_h\in S_h .
\end{align}
For homogeneous equations, i.e. $f\equiv g_j\equiv 0$, the solutions of
(\ref{PDE0}) and (\ref{FEEq0}) can be expressed as $u(t)=E(t)u^0$
and $u_h(t)=E_h(t)u^0_h$, where $\{E(t)=e^{-tA}\}_{t>0}$ and
$\{E_h(t)=e^{-tA_h}\}_{t>0}$ denote the semigroups generated by the
operators $-A$ and $-A_h$, respectively.
From the theory of parabolic equations, we know that
$\{E(t)\}_{t>0}$ extends to an analytic semigroup on
$C(\overline\Omega)$, satisfying
\begin{align}
&\|E(t)v\|_{L^\infty}
+t\|\partial_tE(t)v\|_{L^\infty}\leq C\|v\|_{L^\infty}
,\quad\forall~v\in C(\overline\Omega) .
\end{align}
Its counterpart for the discrete finite element
operator is the analyticity of the semigroup $\{E_h(t)\}_{t>0}$
on $L^\infty\cap S_h$:
\begin{align}
&\|E_h(t)v_h\|_{L^\infty} +t\|\partial_tE_h(t)v_h\|_{L^\infty}\leq
C\|v_h\|_{L^\infty} ,\quad\forall~v_h\in S_h,~ \forall~t>0.
\label{STLEst}
\end{align}
Along with the approach of analytic semigroup, one may reach more
precise analysis of the finite element solution, such as
maximum-norm error estimates of semi-discrete Galerkin FEMs
\cite{Tho, TW1,Wah}, resolvent estimates of elliptic finite element
operators \cite{Bak, BTW,CLT}, error analysis of fully discrete FEMs for
parabolic equations \cite{LN, Pal,Tho}, and the discrete maximal
$L^p$ regularity \cite{Gei1,Gei2}.

A related topic is the space-time maximum-norm
stability estimate for inhomogeneous
equations ($f$ or $g_j$ may not be identically
zero):
\begin{align}
\|u_h\|_{L^\infty(\Omega\times(0,T))}\leq
C_T\|u_h^0\|_{L^\infty}+C_Tl_h\|u\|_{L^\infty(\Omega\times(0,T))},
\quad\forall~T>0
. \label{STLEst2}
\end{align}
Under certain regularity assumptions on $u$, a straightforward
application of the above inequality is the maximum-norm error
estimate:
\begin{align}
\| u - u_h \|_{L^\infty(\Omega\times(0,T))} \leq C_T\| u^0 - u_h^0
\|_{L^\infty} +C_Tl_h h^{r+1}  \| u \|_{L^\infty( (0,T);
W^{r+1,\infty})}  . \label{error}
\end{align}

In the last several decades, many efforts have been devoted to the
stability-analyticity estimate (\ref{STLEst}) and the space-time
stability estimate (\ref{STLEst2}). Schatz et. al. \cite{STW1}
established (\ref{STLEst}) for $d=2$ and $r=1$, with constant
coefficients $a_{ij}$, by using a weighted-norm technique. Later,
Nitsche and Wheeler \cite{NW} proved (\ref{STLEst2}) for $d=2,3$ and
$r\geq 4$ with constant coefficients. Rannacher \cite{Ran} proved
(\ref{STLEst})-(\ref{STLEst2}) in convex polygons with
constant coefficients, and  Chen \cite{Chen} improved the results to
$1\leq d\leq 5$. A more precise analysis was given by Schatz et al
\cite{STW2}, where they proved that \refe{STLEst}-\refe{STLEst2}
hold with $l_h=1$ for $r\geq 2$ and $l_h=\ln(1/h)$ for $r=1$, and
they showed that the logarithmic factor is necessary for $r=1$. In
\cite{STW2}, the proof was given under the condition that the parabolic
Green's function
satisfies 
\begin{align}
|\partial_t^\alpha\partial_x^\beta G(t, x, y)|\leq
C(t^{1/2}+|x-y|)^{-(d+2\alpha+|\beta|)} e^{-\frac{|x-y|^2}{Ct}}
\label{g-cond} ,~~\forall~\alpha\geq 0,|\beta|\geq 0,
\end{align}
which holds when the coefficients $a_{ij}(x)$
are smooth enough \cite{EI70}. The stability
estimate (\ref{STLEst}) was also studied in \cite{Bak, TW2}
for the Dirichlet boundary condition and in \cite{Cro} for
a lumped mass method. Moreover, Leykekhman \cite{Ley}
showed the stability estimate (\ref{STLEst2}) in a more general weighted norm,
and Hansbo \cite{Han} investigated the related $L^s\rightarrow L^r$ stability
estimate.  Also see \cite{TW1,Wah} for some
works in the one-dimensional space. Clearly, all these results were
established for parabolic equations with the coefficient $a_{ij}(x)$
being smooth enough. Related maximum-norm error estimates of
Galerkin FEMs in terms of an elliptic projection and the associated elliptic Green's function
can be found in \cite{BSTW,ELWZ,Lin2, LTW, Nit,Ran,Tho, Whe}.
Some other nonlinear models were analyzed in
\cite{Dob2}. Again, these works were
based on the assumption that the coefficients $a_{ij}$ are smooth
enough.

In many physical applications, the coefficients $a_{ij}$ may
not be smooth enough. One of examples is the
incompressible miscible flow in porous media \cite{Dou,LS1},
where $[a_{ij}]_{i,j=1}^d$ denotes the diffusion-dispersion tensor
which is Lipschitz continuous in many cases. In this case,
the solution is in $L^p((0,T);W^{2,q})$ for $1<p,q<\infty$ 
(see Lemma 2.1 in Section 2).
As a first attempt towards this direction, in this paper,
we prove the maximum-norm stability estimates
(\ref{STLEst})-(\ref{STLEst2}) for parabolic equations with
 Lipschitz continuous coefficients $a_{ij}\in
W^{1,\infty}(\Omega)$, and (\ref{error}) follows immediately. Moreover,
along with these maximum-norm estimates we also obtain a
semigroup estimate:
\begin{align}\label{smgest}
\|\sup_{t>0}|E_h(t)v_h|\big\|_{L^q}\leq C_q\|v_h\|_{L^q},
~~\forall~v_h\in S_h,~~1<q\leq \infty .
\end{align}
Based on these results, we establish the $L^p$ error estimate
\begin{align}
&\|P_hu - u_h \|_{L^p((0,T);L^q)} \leq  C_{p,q}\|P_hu^0 - u_h^0
\|_{L^q} +C_{p,q}\|P_hu-R_hu \|_{L^p((0,T);L^q)} , \label{LpqSt1}
\end{align}
and the maximal $L^p$ regularity
\begin{align}
&\|u_h \|_{L^p((0,T);W^{1,q})} \leq C_{p,q}\|{\bf
g}\|_{L^p((0,T);L^q)} ,\quad\mbox{when}~~ u^0_h\equiv f\equiv 0 ,
\label{LpqSt2}\\
&\|\partial_tu_h \|_{L^p((0,T);L^q)}+\|A_hu_h \|_{L^p((0,T);L^q)}
\leq C_{p,q}\|f\|_{L^p((0,T);L^q)} ,\quad\mbox{when}~~ u^0_h\equiv
{\bf g}\equiv 0, \label{LpqSt3}
\end{align}
for all $1<p,q<\infty$,
where $R_h$ is the Ritz projection operator
associated with the elliptic operator $A$, and $P_h$ is the
$L^2$ projection operator onto the finite element space.
Note that the inequality (\ref{LpqSt3}) was studied by
Geissert \cite{Gei1,Gei2} by assuming that (\ref{g-cond})
holds for $\alpha=0$, $|\beta|\leq 2$ and $0\leq
\alpha\leq 2$, $|\beta|= 2$, where a sufficient condition $a_{ij}\in
C^{2+\alpha}(\overline\Omega)$ was given.
The estimates (\ref{LpqSt2})-(\ref{LpqSt3}) resemble the
maximal $L^p$ regularity of the continuous parabolic
problem. As far as we know, the estimate
(\ref{LpqSt1})-(\ref{LpqSt2}) have not been proved, which imply
optimal error estimates of the finite element solution and
can be regarded as the stability of the parabolic projection
onto the finite element space.
These results are required in \cite{LS2} to establish
optimal $L^p((0,T);L^q)$ error estimates
of FEMs for parabolic equations with time-dependent
nonsmooth coefficients.

The rest part of this paper is organized as follows. In Section 2,
we introduce some notations and present our
main results. In Section 3, we present some new estimates for parabolic
Green's functions. Based on these new estimates, we prove a key lemma
in establishing the the maximum-norm stability estimates.
In Section 4, we prove the maximum-norm stability estimates,
$L^p$ error estimates and maximal $L^p$ regularity for the finite element solution.

\section{Notations and main results}
\setcounter{equation}{0}
Let $\Omega$ be a bounded smooth domain in $\R^d$ ($d\geq 2$).
For any integer $k\geq 0$ and $1\leq p\leq\infty$, let
$W^{k,p}(\Omega)$ be the usual Sobolev space \cite{Adams} of
functions defined in $\Omega$ equipped with the norm
$$
\|f\|_{W^{k,p}(\Omega)}=\left\{
\begin{array}{ll}
\biggl(\sum_{|\beta|\leq k}\int_\Omega|\partial^\beta f|^p\d
x\biggl)^\frac{1}{p},
&1\leq p<\infty
,\\[10pt]
\displaystyle
\sum_{|\beta|\leq k}{\rm ess}\,\sup_{x\in
\Omega}\,|\partial^\beta
f(x)|, & p=\infty,
\end{array}
\right.
$$
where
$$\partial^\beta=
\frac{\partial^{|\beta|}}{\partial
x_1^{\beta_1}\cdots\partial x_d^{\beta_d}}
$$
for the multi-index $\beta=(\beta_1,\cdots,\beta_d)$,
$\beta_1\geq
0$, $\cdots$, $\beta_d\geq 0$, and
$|\beta|=\beta_1+\cdots+\beta_d$.
In particular, we set $L^p(\Omega):=W^{0,p}(\Omega)$ for
$1\leq p\leq\infty$ and
$W^{-k,p}(\Omega):=(W^{k,p'}(\Omega))'$ for any positive
integer $k$, and we set $H^k(\Omega):=W^{k,2}(\Omega)$ for any
integer $k$.

For any integer $k\geq 0$ and $0<\alpha<1$, let
$C^{k+\alpha}(\Omega)$ denote the usual H\"{o}lder
space of functions defined in $\Omega$ equipped
with the norm
$$
\|f\|_{C^{k+\alpha}(\Omega)}=\sum_{|\beta|\leq k}\|D^\beta
f\|_{L^\infty(\Omega)}+\sum_{|\beta|=
k}\sup_{x,y\in\Omega}\frac{|D^\beta f(x)-D^\beta
f(y)|}{|x-y|^\alpha} .
$$

For any Banach space $X$ and a given $T>0$, the Bochner spaces
\cite{Yos} $L^p((0,T);X)$ and $W^{1,p}((0,T);X)$ are equipped with
the norms
\begin{align*}
&\|f\|_{L^p((0,T);X)} =\left\{
\begin{array}{ll}
\displaystyle\biggl(\int_0^T
\|f(t)\|_X^pdt\biggl)^\frac{1}{p},
&
1\leq p<\infty
,\\[10pt]
\displaystyle{\rm ess\,sup}_{t\in(0,T)}\|f(t)\|_X,
& p=\infty,\end{array}
\right.
\\[8pt]
&\|f\|_{W^{1,p}((0,T);X)}=
\|f\|_{L^p((0,T);X)}
+\|\partial_tf\|_{L^p((0,T);X)} ,
\end{align*}
and we set
$ Q_T:=\Omega\times(0,T)$. For
nonnegative integers $k_1$ and $k_2$, we define
\begin{align*}
&\|f\|_{W^{k_1,k_2}_{p,q}( Q_T)}
:=\|f\|_{L^p((0,T);L^{q}(\Omega))}
+\|\partial_t^{k_1}f\|_{L^p((0,T);L^{q}(\Omega))}
+\|f\|_{L^p((0,T);W^{k_2,q}(\Omega))} ,
\end{align*}
and
\begin{align*}
&\|f\|^{(h)}_{W^{k,0}_{p}( Q_T)}
:=\|f\|_{L^p( Q_T)}
+\biggl(\int_0^T\sum_{|\alpha|\leq k}\sum_{l=1}^L
\int_{\tau^h_l}|\partial^\alpha f|^p\d x\d t \biggl)^{\frac{1}{p}} ,
\end{align*}
where $\tau^h_l$, $l=1,\cdots,L$, denote elements of a quasi-uniform triangulation of $\Omega$.

For the simplicity of notations,
in the following sections, we
write $L^p$, $W^{k,p}$, $C^{k+s}$ and
$W^{k_1,k_2}_{p,q}$ as the abbreviations of $L^p(\Omega)$,
$W^{k,p}(\Omega)$, $C^{k+s}(\overline\Omega)$ and
$W^{k_1,k_2}_{p,q}( Q_T)$, respectively. We also set
$L^p( Q_T)=
L^p((0,T);L^p)$,
$W^{k_1,k_2}_{p}
=W^{k_1,k_2}_{p,p}$
for nonnegative integer $k_1,k_2$ and
$1\leq p\leq\infty$, and
 $L^p_h:=L^p\cap S_h$. For any domain $Q\subset
 Q_T$, we define $$Q^t:=\{x\in\Omega:~ (x,t)\in Q\}$$
and
$$
\|f\|_{L^{p,q}(Q)}:=\biggl(\int_0^T\biggl(\int_{Q^t}|f(x,t)|^q\d
x\biggl)^{\frac{p}{q}}\d t\biggl)^{\frac{1}{p}}
,\quad\mbox{for}~ 1\leq p,q<\infty ,
$$
and we use the abbreviations
$$
(\phi,\varphi):=\int_\Omega \phi(x)\varphi(x)\d x,\qquad
[u,v]:=\iint_{ Q_T} u(x,t)v(x,t)\d x\d t .
$$
We write $w(t)=w(\cdot,t)$ as abbreviation for any function $w$ defined on $ Q_T$.

Moreover, we set $a(x)=[a_{ij}(x)]_{d\times d}$ as
a coefficient matrix and define the operators
\begin{align*}
&A:H^1\rightarrow H^{-1},\qquad~\, A_h:S_h\rightarrow S_h,\\
&R_h:H^1\rightarrow S_h, \qquad~~\, P_h:L^2\rightarrow S_h, \\
&\overline\nabla\cdot:(H^1)^d\rightarrow H^{-1},\quad
\overline\nabla_h\cdot:(H^1)^d\rightarrow S_h,
\end{align*}
by
\begin{align*}
&\big(A\phi,v\big)= \big(a\nabla \phi,\nabla
v\big)+\big(c\phi,v\big) \qquad~~~
\mbox{for all~~$\phi ,v \in H^1$},\\[3pt]
&\big(A_h\phi_h,v\big)= \big(a\nabla \phi_h,\nabla
v\big)+\big(c\phi_h,v\big) \quad\,
\mbox{for all~~$\phi_h\in S_h$, $v\in S_h$},\\[3pt]
&\big(A_hR_hw,v\big)=\big(Aw,v\big) \qquad\qquad\qquad~~ \mbox{for
all~~$w \in H^1$ and $v\in S_h$} ,
\\[3pt]
&\big(P_h\phi,v\big)=\big(\phi,v\big) \qquad\qquad\qquad\qquad~~~
\mbox{for all~~$\phi \in L^2$
and $v \in S_h$} ,\\
&\big(\overline\nabla\cdot{\bf w},v\big)=-\big({\bf w},\nabla v\big)
\qquad\qquad\qquad~\, \mbox{for all~~${\bf w} \in (H^1)^d$ and
$v \in H^1$,} \\
&\big(\overline\nabla_h\cdot{\bf w},v\big)=-\big({\bf w},\nabla
v\big) \qquad\qquad\qquad \mbox{for all~~${\bf w} \in (H^1)^d$ and $v
\in S_h$} .
\end{align*}
Clearly,
$R_h$ is the Ritz projection operator associated to the elliptic
operator $A$ and $P_h$ is the $L^2$ projection operator onto the finite element space,
which satisfy
\begin{align*}
&\|u-P_hu\|_{W^{m,p}} \leq C\|u\|_{W^{m,p}},\qquad\quad\, 1\leq p\leq\infty,\\
&\|u-R_hu\|_{W^{m,p}} \leq Ch^{1-m}\|u\|_{W^{1,p}},\quad 1<p\leq\infty,~(m,p)\neq (0,\infty) 
\end{align*}
where $m=0,1$ and $C$ is some positive constant
independent of the mesh size $h$. 

With these notations, for $f\in L^2((0,T);L^2)$, ${\bf g}\in L^2((0,T);L^2)^d$
and $u_0\in L^2$, the problem
\begin{align}\label{sdfjisloudf9890}
\left\{
\begin{array}{ll}
\partial_tu+Au=f-\overline\nabla\cdot{\bf g}
 ,\\[3pt]
u(0)=u_0 ,
\end{array}
\right.
\end{align}
admits a unique solution $u\in L^2((0,T);H^1)\cap
H^1((0,T);H^{-1})$, which coincides with the weak solution of the
initial-boundary
value problem \refe{PDE0}.
The following lemma gives the maximal $L^p$ regularity of the
continuous parabolic problem \cite{KW04}.
\begin{lemma}\label{s00}
{\bf$\!\!\!$(Maximal $\bf L^p$ regularity)}~~\\
{\it Let $u\in L^2((0,T);H^1)\cap H^1((0,T);H^{-1})$
 be the solution to the problem
\begin{align}\label{saf9}
\left\{
\begin{array}{ll}
\partial_tu+Au=f-\overline\nabla\cdot{\bf g},\\[3pt]
u(0)=0   .
\end{array}
\right.
\end{align}
Then the following inequalities hold:
\begin{align}
&\!\!\!\|\partial_tu\|_{L^p((0,T);L^q)}+\|u\|_{L^p((0,T);W^{2,q})}
\leq C_{p,q}\|f\|_{L^p((0,T);L^q)}, ~~~
1<p,q<\infty,~~\mbox{if}~~{\bf g}\equiv 0,
\label{dsk2}\\
&\!\!\!\|u\|_{L^p((0,T);W^{1,q})} \leq C_{p,q}\|{\bf
g}\|_{L^p((0,T);L^q)},\qquad\qquad\qquad\qquad\quad 1<p,q<\infty
,~~\mbox{if}~~f\equiv 0. \label{dsk3}
\end{align}
}
\end{lemma}

Since $L^q((0,T);W^{2,q})\cap 
W^{1,q}((0,T);L^{q})\hookrightarrow
L^p((0,T);W^{1,p})\cap W^{1,p}((0,T);W^{-1,p}) $, 
the inequality \refe{dsk2} implies that
\begin{align}\label{dsk4}
&\|\partial_tu\|_{L^p((0,T);W^{-1,p})}
+\|u\|_{L^p((0,T);W^{1,p})}
\leq C_{p,q}\|f\|_{L^q((0,T);L^q)}, ~~~
\mbox{if}~~{\bf g}\equiv 0,
\end{align} 
for any $p>2$ and $q=(d+2)p/(d+2+p)<p$.

\subsection{Main results}
The main results of this paper are given below and the proofs are presented in
Section 3--4.

\begin{theorem}\label{MainTHM1}
{\it If $a_{ij}=a_{ji}\in W^{1,\infty}$ and $c\in L^\infty$ satisfy the condition {\rm
(\ref{coeffcond})}, then the solutions of
{\rm(\ref{PDE0})}-{\rm(\ref{FEEq0})} satisfy the
maximum-norm estimates {\rm(\ref{STLEst})}-{\rm(\ref{STLEst2})} with $l_h=\ln(2+1/h)$,
the semigroup estimate {\rm(\ref{smgest})}, the $L^p$ stability
estimate {\rm(\ref{LpqSt1})} and the maximal $L^p$
regularity {\rm(\ref{LpqSt2})}-{\rm(\ref{LpqSt3})}.
}
\end{theorem}

In the proof of our main results, we can assume that the functions $f$ and $g$ are smooth enough and the exact solution $u$ satisfies $u\in L^p((0,T);W^{2,p})$ and $\partial_tu\in L^p((0,T);L^p)$ for arbitrarily large $p$. However, the generic positive constant $C$ in this paper does not depend on the regularity of $f$, $g$ or $u$.
Therefore, by a passing to a limit, one can see that
(1.8) defines a parabolic projection for $u^0_h\in S_h$ and 
$u\in C(\overline \Omega\times[0,T])$,
(1.13) holds for ${\bf g}\in L^p((0,T);L^q)^d$ and (1.14)
holds for $f\in L^p((0,T);L^q)$.

Unlike \cite{STW2},
we are not able to remove the logarithmic factor $l_h$ in 
{\rm(\ref{STLEst})}-{\rm(\ref{STLEst2})} for finite element spaces 
of polynomial 
degree $r\geq 2$, due to the low regularity of the coefficients.

\subsection{Further notations}\label{fnot}
To prove our main results, we present some further notations,
which were introduced in \cite{STW2,SW1} and also used in 
\cite{Ley,Sol}.

For an element $\tau_l^h$ and a point $x_0\in \tau_l^h$, we
let $\delta_{x_0}$ denote the Dirac Delta function centered at
$x_0$, i.e. $\int_\Omega\delta_{x_0}(y)\varphi (y)\d y=\varphi(x_0)$ for
any $\varphi\in C(\overline\Omega)$, and we denote by
$\widetilde\delta_{x_0}$ a regularized
Delta function satisfying the following
conditions:
\begin{align}
&\mbox{$\widetilde\delta_{x_0}$ is
supported in $\tau_l^h$ ,}\\
&\chi(x_0)=\int_{\tau_l^h}\chi\,
\widetilde\delta_{x_0}\,\d x,\quad
\mbox{for all}~\chi\in S_h ,
\label{chidef}\\
&\|\widetilde\delta_{x_0}\|_{W^{m,p}}\leq
 Ch^{-m-d(1-1/p)}~~\mbox{for}~~1\leq p\leq
 \infty,~m=0,1,2,3 . \label{Deltadef}
\end{align}
Let $G(t,x,x_0)$ be {\it Green's function} of the parabolic
equation, defined by
\begin{align}\label{GFdef}
&\partial_tG(t,\cdot,x_0)+ AG(t,\cdot,x_0)=0\quad
\mbox{for $t>0$ with $G(0,\cdot,x_0)=\delta_{x_0}$}
.
\end{align}
The corresponding {\it regularized Green's function} $\Gamma(t,x,x_0)$
is defined by
\begin{align}\label{GMFdef}
&\partial_t\Gamma(\cdot,\cdot,x_0)+A\Gamma(\cdot,\cdot,x_0)=0\quad
\mbox{for~ $t>0$~ with~
$\Gamma(0,\cdot,x_0)=\widetilde\delta_{x_0}$},
\end{align}
and the {\it discrete Green's function} $\Gamma_h(\cdot,\cdot,x_0)$
is defined as the solution of the equation
\begin{align}
&\partial_t\Gamma_{h}(\cdot,\cdot,x_0)
+ A_h\Gamma_{h}(\cdot,\cdot,x_0)=0 \label{GMhFdef}
\quad\mbox{for~ $t>0$~ with~
$\Gamma_{h}(0,\cdot,x_0)=P_h\delta_{x_0}=P_h\widetilde\delta_{x_0}$},
\end{align}
where $P_h$ is the $L^2$ projection onto the finite element space.
Note that
$\Gamma(t,x,x_0)$ and $\Gamma_h(t,x,x_0)$ are symmetric with respect
to $x$ and $x_0$.

By the fundamental estimates of parabolic equations \cite{FS}
and from Appendix B of \cite{Gei1}, we know that the
Green's function $G$ satisfies
\begin{align}
&|G(t,x,y)|\leq
C(t^{1/2}+|x-y|)^{-d}e^{-\frac{|x-y|^2}{Ct}},\label{FEstP}\\
&|\partial_tG(t,x,y)|\leq Ct^{-d/2-1}e^{-\frac{|x-y|^2}{Ct}} ,\label{FtEstP}\\
&|\partial_{tt}G(t,x,y)|\leq Ct^{-d/2-2}e^{-\frac{|x-y|^2}{Ct}} .\label{FtEstP2}
\end{align}
By estimating $\Gamma(t,x,x_0)=\int_\Omega
G(t,x,y)\widetilde\delta_{x_0}(y)\d y$, we see that
(\ref{FEstP})-(\ref{FtEstP2}) also hold when $G$ is replaced by
$\Gamma$ and when $\max(t^{1/2},|x-y|)\geq 2h$.

For any open subset $D\subset\Omega$, we set $\overline
D^\partial ={\rm int}(D)\cup (\overline D \cap
\partial\Omega)$. Let $S_h(D)$ denote the restriction of
functions in $S_h$ to $D$, and let $S_h^0(\overline
D^\partial)$ denote the functions in $S_h$ with the support in
$\overline D^\partial$. For a given subset $D\subset\Omega$,
we set $D_{d} =\{x\in \Omega : {\rm dist}(x,D)\leq d\}$
for $d> 0$.
We denote by $I_h : W^{1,1}(\Omega) \rightarrow S_h$ the
operator given in \cite{STW2} having the following properties:
if $d \geq kh$, then
$$\| I_hv - v\|_{W^{s,p}(D_d)}^{(h)} \leq
Ch^{l-s}\|v\|_{W^{l,p}(D_{2d})},
\quad~\mbox{for~ $0\leq s\leq l\leq r $ ~and~ $ 1\leq p\leq
\infty$,}$$
and if supp$(v)\subset \overline D ^\partial_d$, then
$I_hv\in S^0_h(\overline D_{2d}^\partial)$; also, if
$v|_{D_d}\in S_h(D_d)$, then $I_hv = v $ on $D$ and the
bound above may be replaced by
$Ch^{l-s}\|v\|_{W^{l,p}(D_{2d}\backslash D)}$.

For any integer $j$, we define $d_j=2^{-j}$.
Let $J_1=1$ and $J_0=0$, and let $J_*$ be an integer
satisfying $2^{-J_*}= C_*h$ with $C_*\geq 16$ to
be determined later, thus $J_*=\log_2[1/(C_*h)]\leq 2\ln(2+1/h)$.
For the given constant $C_*$, we have $J_1<J_*$ when $h<1/(4C_*)$,
and for a given $x_0\in\Omega$ and $j\geq
J_1$ we define the subsets $Q_*,Q_j\subset\Omega_T$ and
$\Omega_*,\Omega_j\subset\Omega$ by
\begin{align*}
&Q_*(x_0)=\{(x,t)\in\Omega_T: \max (|x-x_0|,t^{1/2})\leq
d_{J_*}\},\\
&\Omega_*(x_0)=\{x\in \Omega: |x-x_0|\leq d_{J_*}\},\\
&Q_j(x_0)=\{(x,t)\in \Omega_T:d_j\leq
\max (|x-x_0|,t^{1/2})\leq2d_j\},\\
&\Omega_j(x_0)=\{x\in \Omega: d_j\leq|x-x_0|\leq2d_j\} ,\\
&Q_{J_0}(x_0)=\Omega_T\big\backslash\big(
\cup_{j=J_1}^{J_*}Q_{j}(x_0)\cup Q_*(x_0)\big) ,\\
&\Omega_{J_0}(x_0)=\Omega\big\backslash\big(
\cup_{j=J_1}^{J_*}\Omega_{j}(x_0)\cup \Omega_*(x_0)\big).
\end{align*}
For $j<J_0$, we simply define
$Q_{j}(x_0)=\Omega_{j}(x_0)=\emptyset$ and for any integer $j$ we define
\begin{align*}
&\Omega_j'(x_0)=\Omega_{j-1}(x_0)\cup\Omega_{j}(x_0)\cup\Omega_{j+1}(x_0),
~~~Q_j'(x_0)=Q_{j-1}(x_0)\cup Q_{j}(x_0)\cup Q_{j+1}(x_0), \\
&\Omega_j''(x_0)=\Omega_{j-2}(x_0)\cup\Omega_{j}'(x_0)\cup\Omega_{j+2}(x_0),
~~~
Q_j''(x_0)=Q_{j-2}(x_0)\cup Q_{j}'(x_0)\cup Q_{j+2}(x_0),\\
&\Omega_j'''(x_0)=\Omega_{j-2}(x_0)\cup\Omega_{j}''(x_0)\cup\Omega_{j+2}(x_0),
~~~
Q_j'''(x_0)=Q_{j-2}(x_0)\cup Q_{j}''(x_0)\cup Q_{j+2}(x_0) .
\end{align*}
Then we have
\begin{align*}
&\Omega_T=\bigcup^{J_*}_{j=J_0}Q_j(x_0)\,\cup Q_*(x_0)
\quad\mbox{and}\quad
\Omega=\bigcup^{J_*}_{j=J_0}\Omega_j(x_0)\,\cup \Omega_*(x_0),\end{align*}
We refer to $Q_*(x_0)$ as the ``innermost" set.
We shall write $\sum_{*,j}$ when the innermost set is included and
$\sum_j$ when it is not. When $x_0$ is fixed and there is no ambiguity, we simply write
$Q_j=Q_j(x_0)$, $Q_j'=Q_j'(x_0)$, $Q_j''=Q_j''(x_0)$ and
$\Omega_j=\Omega_j(x_0)$, $\Omega_j'=\Omega_j'(x_0)$,
$\Omega_j''=\Omega_j''(x_0)$.

We shall write
\begin{align*}
&\|v\|_{k,D}=\biggl(\int_D\sum_{|\alpha|\leq
k}|\partial^\alpha v|^2\d x\biggl)^{\frac{1}{2}},\qquad
\vertiii{v}_{k,Q}=\biggl(\int_Q\sum_{|\alpha|\leq
k}|\partial^\alpha v|^2\d x\d t\biggl)^{\frac{1}{2}} ,
\end{align*}
for any domain $D\subset\Omega$ and $Q\subset\Omega\times(0,T)$.
The time derivative will always be displayed explicitly.
We denote by $C$ a generic positive constant, which
will be independent of $h$, $x_0$, and
the undetermined constant $C_*$ until
it is determined at the end of Section \ref{dka7}.

\section{Estimates of the parabolic Green's function}
\setcounter{equation}{0}
To prove our main results, we need the following lemma.
The proof of the lemma will be given in the next two subsetions.
\begin{lemma}\label{GMhEst}
{\it Let $x_0\in\Omega$ and $T=1$.
Let $\Gamma(t,x,x_0)$ and $\Gamma_h(t,x,x_0)$ be defined in
{\rm(\ref{GMFdef})-(\ref{GMhFdef})}, and set
$F(t,x)=\Gamma_h(t,x,x_0)-\Gamma(t,x,x_0)$.
Then there exists a positive constant $h_0>0$
such that when $h<h_0$ we have
\begin{align}
&
\|\partial_tF\|_{L^1( Q_T)}+
\|t\partial_{tt}F\|_{L^1( Q_T)}
+
h^{-1}l_h^{-1}
\|F\|_{W^{1,0}_1( Q_T)}\leq C,\label{FFEst1}\\
& \|\partial_tF\|_{L^1(\Omega\times(0,\infty))}
+ \|t\partial_{tt}F\|_{L^1 (\Omega\times(0,\infty))} \leq C ,
\label{FFEst2}
\end{align}
where $l_h=\ln(2+1/h)$ and the constant $C$ does not depend on $x_0$.
}
\end{lemma}

The estimates in the lemma were proved in \cite{STW2} 
for parabolic equations with smooth coefficients 
for which the Green function satisfies (\ref{g-cond}).
Since $x_0$ is fixed, we simply write
$G$ and $\Gamma$ as abbreviations for
the functions $G(\cdot,\cdot,x_0)$ and
$\Gamma(\cdot,\cdot,x_0)$, respectively,
when there is no ambiguity.  We shall assume that
$h< 1/(4C_*)$, so that $Q_j(x_0)$, $J_0\leq j\leq J_*$,
are well defined as in the last section.
In the rest part of this section, we set $T=1$.

\subsection{Estimates of the Green's functions}
\label{EGF}
In this subsection, we present some new estimates for the Green's
function, the regularized Green's function and the discrete Green's
function, which will be used in the next subsection to prove Lemma \ref{GMhEst}.

\begin{lemma}\label{GFEst1}
{\it
There exist $p_1>d$ and $\alpha>0$ such that 
for any integer $J_0\leq j\leq J_*$, we have
\begin{align}
&d_j^{-4}\vertiii{\Gamma(\cdot,\cdot,x_0)}_{2,Q_j(x_0)}
+d_j^{-2}\vertiii{\partial_{t}
\Gamma(\cdot,\cdot,x_0)}_{2,Q_j(x_0)} \nn\\
&\qquad\qquad\qquad\qquad\qquad\qquad
+\vertiii{\partial_{tt}\Gamma(\cdot,\cdot,x_0)}_{2,Q_j(x_0)}
\leq Cd_j^{-d/2-5}, \label{GFest01}\\
&\|\partial_t\partial_{x_i}G(\cdot,\cdot,x_0)
\|_{L^\infty(Q_j(x_0))}+\|\partial_t
\partial_{x_i}\Gamma(\cdot,\cdot,x_0)\|_{L^\infty(Q_j(x_0))}\leq
Cd_j^{-d -3},\label{GFest02}\\
&\|\partial_{x_i x_l}G(\cdot,\cdot,x_0) \|_{L^{\infty,p_1}(\cup_{k\leq
j}Q_k(x_0))}\leq Cd_j^{-d-2+d/p_1} \label{GFest03}
,\\
&\|\partial_t\Gamma(\cdot,\cdot,x_0)\|_{L^1(\Omega\times(T,\infty))}
+\|t\partial_{tt}\Gamma(\cdot,\cdot,x_0)\|_{L^1(\Omega\times(T,\infty))}\leq
C ,\label{GFest0423}\\
&d_j^{-\alpha}\|\partial_{x_i y_l}
G(\cdot,\cdot,x_0)\|_{L^\infty(Q_k(x_0))} +\|\partial_{x_i
y_l}G(\cdot,\cdot,x_0)\|_{C^{\alpha,\alpha/2}(\overline
Q_k(x_0))}\leq Cd_j^{-d -2-\alpha}\label{GFest05} ,\\
&\|\Gamma_h(1,\cdot,x_0)\|_{L^2}+\|\partial_t\Gamma_h(1,\cdot,x_0)\|_{L^2}
+\|\partial_{tt}\Gamma_h(1,\cdot,x_0)\|_{L^2}\leq
C\|\Gamma_h(\cdot,\cdot,x_0)\|_{L^2(\Omega\times(1/2,1])} ,
\label{dlky60}
\end{align}
for $i,l=1,2,\cdots,d$.
}
\end{lemma}
\noindent{\it Proof}~~~
For the given $x_0$ and $j$,  we define a coordinate transformation
$x-x_0=d_j\widetilde x$ and $t=d_j^2\widetilde t$, and $\widetilde
G(\widetilde t,\widetilde x):=G(t,x,x_0)$, $\widetilde
G_{y_l} (\widetilde t,\widetilde x):=\partial_{y_l} G(t,x,y)|_{y=x_0}$, $\widetilde
a(\widetilde x):=a(x)$, 
$\widetilde
c(\widetilde x):=c(x)$,
$\widetilde Q_k=\{(\widetilde x,\widetilde
t)\in\R^{d+1}:(x,t)\in Q_k\}$, $\widetilde Q_k'=\widetilde
Q_{k-1}\cup \widetilde Q_k\cup \widetilde Q_{k+1}$,  
$\widetilde \Omega_k=\{(\widetilde x,\widetilde
t)\in\R^{d+1}:(x,t)\in \Omega_k\}$, $\widetilde \Omega_k'=\widetilde
\Omega_{k-1}\cup \widetilde \Omega_k\cup \widetilde \Omega_{k+1}$,
$\widetilde \Omega=\{\widetilde x\in\R^d: x \in
\Omega\}$, and $\widetilde Q_{\widetilde T}=\{(\widetilde
x,\widetilde t)\in\R^{d+1}: (x,t) \in  Q_T\}$.
Then $\widetilde G(\widetilde t,\widetilde
x)$ and $\widetilde G_{y_l} (\widetilde t,\widetilde
x)$ are solutions of the equations
\begin{align}
&\partial_{\widetilde t}\widetilde G-\nabla_{\widetilde
x}\cdot(\widetilde a\nabla_{\widetilde x}
\widetilde G)+\widetilde c\widetilde G=0 \quad\mbox{in}~~\widetilde Q_{j}', \label{GFest06}\\
&\partial_{\widetilde t}\widetilde G_{y_l}-\nabla_{\widetilde
x}\cdot(\widetilde a\nabla_{\widetilde x}
\widetilde G_{y_l})
+\widetilde c\widetilde G_{y_l}=0 \quad\mbox{in}~~\widetilde Q_{j}' .
\label{GFest066}
\end{align}

By the estimates of parabolic equations 
(see Lemma A.1 in Appendix), 
we have
\begin{align}
&|\!|\!|\partial_{\widetilde t}\widetilde
G|\!|\!|_{\widetilde Q_j}+|\!|\!|\widetilde
G|\!|\!|_{2,\widetilde Q_j}
+|\!|\!|\partial_{\widetilde t}\widetilde
G|\!|\!|_{2,\widetilde Q_j}
+|\!|\!|\partial_{\widetilde t\widetilde t}\widetilde
G|\!|\!|_{2,\widetilde Q_j}
\leq  C|\!|\!|\widetilde G|\!|\!|_{\widetilde Q_j'}  , \label{dkjq}\\
&\|\partial_{\widetilde x_i}\widetilde
G\|_{L^\infty(\widetilde Q_j)}
+\|\partial_{\widetilde x_i}\widetilde
G\|_{C^{\alpha,\alpha/2}(\overline{\widetilde Q}_j)}
+\|\partial_{\widetilde x_i}\partial_{\widetilde x_l}\widetilde
G\|_{L^{\infty,p_1}(\widetilde Q_j)}
 \leq C|\!|\!|\widetilde G|\!|\!|_{\widetilde Q_j'}
\\
&\|\partial_{\widetilde x_i}\widetilde
G_{y_l}\|_{L^\infty(\widetilde Q_j)} 
+\|\partial_{\widetilde x_i}\widetilde
G_{y_l}\|_{C^{\alpha,\alpha/2}(\overline{\widetilde Q}_j)}
\leq C|\!|\!|\widetilde G_{y_l}|\!|\!|_{\widetilde Q_j'}  .
\label{dkjq3}
\end{align}
 
Transforming back to the $(x,t)$ coordinates, \refe{dkjq}-\refe{dkjq3}
reduce to
\begin{align*}
&d_j^{-4}\vertiii{ G}_{2,Q_j}
+d_j^{-2}\vertiii{\partial_{t}
G}_{2,Q_j}+\vertiii{\partial_{tt}G}_{2,Q_j}\leq
Cd_j^{-6}\vertiii{ G}_{Q_j'} ,\\
& d_j^{-2}\|\partial_{x_i}G\|_{L^\infty( Q_j)}
+\|\partial_t\partial_{x_i}G\|_{L^\infty( Q_j)}\leq Cd_j^{-d/2
-4}\vertiii{ G}_{Q_j'},\\
&
\|\partial_{x_i}\partial_{x_l}G\|_{L^{\infty,p_1}( Q_j)}\leq Cd_j^{-d/2
-3+d/p_1}\vertiii{ G }_{Q_j'} ,\\
& d_j^{-\alpha}\|\partial_{x_i}
\partial_{y_l}G\|_{L^\infty(Q_j)}
+\|\partial_{x_i}\partial_{y_l}G
\|_{C^{\alpha,\alpha/2}(\overline
Q_j)} \leq Cd_j^{-d/2 -2-\alpha}
\vertiii{ \partial_{y_l}G }_{Q_j'}\leq Cd_j^{-1-\alpha}
\|\partial_{y_l}G \|_{L^\infty(Q_j')}.
\end{align*}
From the Green function estimate (\ref{FEstP}), we see that
$\vertiii{ G}_{Q_j'}\leq Cd_j^{-d/2+1}$ and so
\begin{align}
&d_j^{-4}\vertiii{ G}_{2,Q_j}
+d_j^{-2}\vertiii{\partial_{t}
G}_{2,Q_j}+\vertiii{\partial_{tt}G}_{2,Q_j}
\leq Cd_j^{-d/2-5}, \label{GFest07}\\
&d_j^{-2}\|\partial_{x_i}G\|_{L^\infty(Q_j)}
+\|\partial_t\partial_{x_i}G\|_{L^\infty(Q_j)}\leq Cd_j^{-d
-3},\label{GFes03}\\
&\|\partial_{x_i x_l}G\|_{L^{\infty,p_1}(Q_j)}\leq Cd_j^{-d
-2+d/p_1} \label{Gsnn} ,\\
&d_j^{-\alpha}\|\partial_{x_i}\partial_{y_l}
G\|_{L^\infty(Q_j)}+\|\partial_{x_i}\partial_{y_l}
G\|_{C^{\alpha,\alpha/2}(\overline Q_j)} \leq Cd_j^{
-1-\alpha}\|\partial_{y_l}G \|_{L^\infty(Q_j')}
\leq
Cd_j^{-d-2-\alpha} , \label{Gmm2}
\end{align}
where we have used (\ref{GFes03}) in deriving (\ref{Gmm2}).
Clearly, (\ref{Gsnn}) further implies that
\begin{align}\label{fd00}
&\|\partial_{x_i x_l}G\|_{L^{\infty,p_1}(\cup_{k\leq j}Q_k)}\leq
Cd_j^{-d -2+d/p_1} .
\end{align}

By estimating $\Gamma(t,x)=\int_\Omega
G(t,x,y)\widetilde\delta_{x_0}(y)\d y$, we can see that the
estimates (\ref{GFest07})-(\ref{fd00}) also hold when $G$ is
replaced by $\Gamma$.

From the inequalities
(\ref{FtEstP})-(\ref{FtEstP2}) we derive that
\begin{align}\label{dlk80}
&\|\partial_tG(t,\cdot,x_0)\|_{L^\infty}
+\|\partial_t\Gamma (t,\cdot,x_0)\|_{L^\infty}
\nn\\
&+\|t\partial_{tt}G(t,\cdot,x_0)\|_{L^\infty}
+\|t\partial_{tt}\Gamma (t,\cdot,x_0)\|_{L^\infty}\leq
Ct^{-d/2-1}~~\mbox{for}~~t\geq 1/4,
\end{align}
which implies (\ref{GFest0423}).

Finally, we note that the inequality (\ref{dlky60}) follows
from basic energy estimates.

The proof of Lemma \ref{GFEst1} is complete.
~\endproof\bigskip

\subsection{Proof of Lemma \ref{GMhEst}}
\label{dka7}

\begin{lemma}\label{LocEst}
{\it
Suppose that $z(\cdot,t)\in H^1$, $z_t(\cdot,t)\in L^2$ and $z_h(\cdot,t)\in S_h$
for each fixed $t\in[0,T]$, and suppose that $e=z_h-z$ satisfies the equation
$$
(e_t,\chi)+(a\nabla e,\nabla \chi)+(c
e,\chi)=0,\quad\forall~\chi\in S_h,~t>0 ,
$$
with $z(\cdot,0)=0$ and $z_h(\cdot,0)=z_{0h}$ on $\Omega_j'$.
Then for any $q>0$ there exists a constant $C_q$ such that
\begin{align*}
\vertiii{e_t}_{Q_j} + d_j^{-1}\vertiii{e}_{1,Q_j}
\leq
C_q\big(I_j(z_{0h})+X_j(I_hz-z)
+H_j(e)+d_j^{-2}\vertiii{e}_{Q_j'}\big),
\end{align*}
where
\begin{align*}
&I_j(z_{0h})=\|z_{0h}\|_{1,\Omega_j'}
+d_j^{-1}\|z_{0h}\|_{\Omega_j'},\\
&X_j(I_hz-z)=d_j\vertiii{\partial_t(I_hz-z)}_{1,Q_j'}
+\vertiii{\partial_t(I_hz-z)}_{Q_j'}
+d_j^{-1}\vertiii{I_hz-z}_{1,Q_j'}+
d_j^{-2}\vertiii{I_hz-z}_{Q_j'},\\
&H_j(e)=(h/d_j)^q\big(\vertiii{e_t}_{Q_j'}
+d_j^{-1}\vertiii{e}_{1,Q_j'}\big) .
\end{align*}
}
\end{lemma}

The above lemma was proved in \cite{STW2} (Section 5 and Section 6) only for
parabolic equations with smooth coefficients. However, we can see from the proof
that the lemma still holds when $a_{ij}\in W^{1,\infty}(\Omega)$ and
$c\in L^\infty(\Omega)$ satisfy (\ref{coeffcond}).
Moreover, for parabolic equations with smooth coefficients,
Lemma \ref{GMhEst} was proved in \cite{STW2} by applying Lemma
\ref{LocEst} with the additional assumption (\ref{g-cond}). Here, we
shall prove Lemma \ref{GMhEst} directly from Lemma \ref{GFEst1} and
Lemma \ref{LocEst}.

First, we prove \refe{FFEst1}. 
Let $\mu_j=[h\ln(2+1/h)]^{-1}+d_j^{-1}$ and define
\begin{align}\label{KdKj0}
K_j=
\vertiii{\partial_tF}_{Q_j}+
d_j^2\vertiii{\partial_{tt}F}_{Q_j}
+
\mu_j
\vertiii{F}_{1,Q_j} ,
\end{align}
and
\begin{align}\label{KdKj}
{\cal K}:=\sum_{j}d_j^{1+d/2}K_j.
\end{align}
From Section 4 of \cite{STW2} we see that (\ref{FFEst1}) holds if we can prove that ${\cal K}\leq C$ for some positive constant $C$ which is independent of $h$, $J_*$ and $C_*$.

To prove the boundedness of ${\cal K}$, we set $e=F$ and $e=F_t$
in Lemma \ref{LocEst}, respectively. Since in either case $z(0)=0$ on $\Omega_j'$,
we obtain
\begin{align}
K_j\leq
C(\widehat{I_j}+\widehat{X_j}+\widehat{H_j}+\mu_jd^{-1}_j\vertiii{
F }_{Q'_j} ),
\end{align}
where, by using the exponential decay of 
$|P_h\widetilde\delta_{x_0}(y)|\leq Ch^{-d}e^{-C|y-x_0|/h}$ 
\cite{STW2}, we have
\begin{align*}
&\widehat{I_j}=d_j^2\|F_t(0)\|_{1,\Omega_j'}
+d_j\|F_t(0)\|_{\Omega_j'}
+d_j\mu_j\|F(0)\|_{1,\Omega_j'}+\mu_j\|F(0)\|_{\Omega_j'}
\leq Ch^{-1-d/2}e^{-Cd_j/h},\\
&\widehat{X_j}=
d_j^3\vertiii{(I_h\Gamma-\Gamma)_{tt}}_{1,Q_j'}
+d_j^2\vertiii{(I_h\Gamma-\Gamma)_{tt}}_{Q_j'}
+\mu_jd_j^2\vertiii{(I_h\Gamma-\Gamma)_t}_{1,Q_j'}
+\mu_jd_j\vertiii{(I_h\Gamma-\Gamma)_t}_{Q_j'}\\
&\qquad~
+\mu_j\vertiii{I_h\Gamma-\Gamma}_{1,Q_j'}
+\mu_jd_j^{-1}\vertiii{I_h\Gamma-\Gamma}_{Q_j'}\\
&\quad~
\leq (d_j^3h+d_j^2h^2)\vertiii{\partial_{tt}\Gamma}_{2,Q_j''}
+(\mu_jd_j^2h+\mu_jd_jh^2)\vertiii{\partial_{t}\Gamma}_{2,Q_j''}
+(\mu_jh+\mu_jd_j^{-1}h^2)\vertiii{\Gamma}_{2,Q_j''}\\
&\quad~
\leq C hd_j^{-d/2-2}+C(
[\ln(2+1/h)]^{-1}+h/d_j)d_j^{-d/2-1},\\
&\widehat{H_j}=\big(h/d_j\big)^q \big(d_j^2\vertiii{F_{tt}}_{Q_j'}
+d_j\vertiii{F_{t}}_{1,Q_j'}
+\mu_jd_j\vertiii{F_t}_{Q_j'}+\mu_j\vertiii{F}_{1,Q_j'}
\big)\\
&\quad~
\leq \big(h/d_j\big)^q
\big(d_j^2\vertiii{F_{tt}}_{Q_T}+d_j\vertiii{F_{t}}_{1,Q_T}
+\mu_jd_j\vertiii{F_t}_{Q_T}+\mu_j\vertiii{F}_{1,Q_T}
\big).
\end{align*}
The last term $\widehat{H_j}$ was estimated in \cite{STW2} via
energy estimates, with $\sum_{j}d_j^{1+d/2}\widehat{H_j}\leq C$.
Therefore,
\begin{align}\label{dlj6}
{\cal K}=\sum_{j}d_j^{1+d/2}K_j\leq
C+C\sum_{j}d_j^{d/2}\mu_j\vertiii{F}_{Q_j'} .
\end{align}

To estimate $\vertiii{F}_{Q_j'}$, we apply a duality argument.
Let $w$ be the solution of the backward parabolic equation
$$
-\partial_tw+Aw=v\quad\mbox{with}~~w(T)=0,
$$
where $v$ is a function which is supported in $Q_j'$ and
$\vertiii{v}_{Q_j'}=1$. Multiplying the above equation by $F$,
with integration by parts we get
\begin{align}\label{dka6}
[F,v]=(F(0),w(0))+[F_t,w]+\sum_{i,j=1}^d[a_{ij}\partial_j
F,\partial_i w]+[cF,w],
\end{align}
where
\begin{align*}
(F(0),w(0))&=(P_h\widetilde\delta_{x_0}-\widetilde\delta_{x_0},w(0))\\
&=(P_h\widetilde\delta_{x_0}-\widetilde\delta_{x_0},w(0)-I_hw(0))\\
&=
(P_h\widetilde\delta_{x_0},w(0)-I_hw(0))_{\Omega_j''}
+(P_h\widetilde\delta_{x_0}-\widetilde\delta_{x_0},
w(0)-I_hw(0))_{(\Omega_j'')^c}\\
&:=I_1+I_2 .
\end{align*}
Since
$|P_h\widetilde\delta_{x_0}(y)|\leq Ch^{-d}e^{-C|y-x_0|/h}$ 
\cite{STW2}, we derive
that
\begin{align}
&|I_1|\leq
Ch\|P_h\widetilde\delta_{x_0}\|_{L^2(\Omega_j'')}\|w(0)\|_{H^1(\Omega)}\leq
Cd_j^{d/2}h^{-d+1}e^{-Cd_j/h}\vertiii{v}_{Q_j'}\leq
Ch^{2}d_j^{-d/2 -1}, \label{SF2}\\
&|I_2|\leq C\|P_h\widetilde\delta_{x_0}-
\widetilde\delta_{x_0}\|_{L^{p_1'}}
\|w(0)-I_hw(0)\|_{L^{p_1}((\Omega_j'')^c)}\leq
Ch^{2-d/p_1}\|w(0)\|_{W^{2,p_1}((\Omega_j'')^c)} .
\label{SF22}
\end{align}
We proceed to estimate $\|w(0)\|_{W^{2,p_1}((\Omega_j'')^c)}$.
Let $D_j$ be a set containing $(\Omega_j'')^c$ but its
distance to $\Omega_j'$ is larger than $C^{-1}d_j$. Since
$$
\partial_{x_i}\partial_{x_j}w(x,0)
=\int_0^{T}\int_{\Omega}
\partial_{x_i}\partial_{x_j}G(s,x,y)v(y,s)\d y\d s ,
$$
by taking the $L^p(D_j)$ norm with respect to $x$ we obtain
$$
|x-y|+s^{1/2}\geq C_1^{-1}d_j \quad\mbox{for $x\in D_j$ and $(y,s)\in
Q_j$}
$$
for some positive constant $C_1$.  Using (\ref{GFest03}) we
further derive that
\begin{align}
\|\partial_{x_i}\partial_{x_j}w(0)\|_{L^{p_1}(D_j)}
&\leq C\sup_{y\in\Omega}\|\partial_{x_i}\partial_{x_j}
G(\cdot,\cdot,y)\|_{L^{\infty,p_1}(\cup_{k\leq
j+\log_2C_1}Q_k(y))}\|v\|_{L^{1}(Q_j')}\nn\\
&\leq C d_j^{-d -2+d/p_1}\|v\|_{L^{1}(Q_j')}
 \nn\\
&\leq C d_j^{-d/2 -1+d/p_1}\vertiii{v}_{Q_j'} , \nn\\
&=C d_j^{-d/2 -1+d/p_1} . \label{SF3}
\end{align}
From (\ref{SF2})-(\ref{SF3}), we see that the 
first term on the right-hand side of \refe{dka6} is bounded by
\begin{align}
|(F(0),w(0))| \leq  Ch^{2}d_j^{-d/2 -1}+Ch^{2}d_j^{-d/2
-1}(h/d_j)^{-d/p_1} \leq  Ch^{2}d_j^{-d/2 -1}(h/d_j)^{-d/p_1}  ,
\label{f69}
\end{align}
and the rest terms are bounded by
\begin{align}\label{sd80}
&[F_t,w]+\sum_{i,j=1}^d[a_{ij}\partial_j F,\partial_iw]+[cF,w]\nn\\
&=[F_t,w-I_hw]+\sum_{i,j=1}^d[a_{ij}\partial_j F,\partial_i
(w-I_hw)] +[cF,w-I_hw]\nn\\
&\leq
\sum_{*,i}C(h^2\vertiii{F_t}_{Q_i}+h\vertiii{F}_{1,Q_i})\vertiii{w}_{2,Q_i'}
.
\end{align}
To estimate $\vertiii{w}_{2,Q_i'}$ we consider the expression
$$
\partial_{x_i}\partial_{x_j}w(x,t)
=\int_0^{T}\int_{\Omega}
\partial_{x_i}\partial_{x_j}G(s-t,x,y)v(y,s)1_{s>t}\,\d y\d s .
$$
If $i\leq j-2$ (so that $d_i>d_j$), then $w(t)=0$ for
$t>d_j^2$ (because $v$ is supported in $Q_j$), $|x-y|\sim d_i$
and $s-t\in(0,d_i^2)$ for $t< d_j^2$, $(x,t)\in Q_i$ and
$(y,s)\in Q_j'$.
we obtain
\begin{align*}
\vertiii{\partial_{x_i}\partial_{x_j}w}_{Q_i'}
\leq\sup_{y}\vertiii{\partial_{x_i}\partial_{x_j}G(\cdot,\cdot,y)}_{Q_i(y)}
\|v\|_{L^1(Q_j)}\leq
Cd_i^{-d/2-1}d_j^{d/2+1}\vertiii{v}_{Q_j}\leq
C(d_j/d_i)^{d/2+1} .
\end{align*}
If $i\geq j+2$ (so that $d_i\leq d_j$), then
$\max(|s-t|^{1/2},|x-y|)\geq d_{j+2}$ for $(x,t)\in Q_i$, thus for
$1/2=1/\bar p_1+1/p_1$ we have
\begin{align*}
\vertiii{\partial_{x_i}\partial_{x_j}w}_{Q_i'}
&=\sup_{(y,s)\in Q_T}\vertiii{\partial_{x_i}\partial_{x_j}G(s-\cdot,\cdot,y)
1_{\cup_{k\leq j+2}Q_{k}(y)}}_{Q_i'}
\|v\|_{L^1(Q_j')}\\
&\leq Cd_i^{ 1+d /\bar p_1}\sup_{y}\|\partial_{x_i}
\partial_{x_j}G(\cdot,\cdot,y)\|_{L^{\infty,p_1}(\cup_{k\leq j+2}Q_{k}'(y))}
\vertiii{v}_{Q_j'}d_j^{d/2+1}
\\
&\leq Cd_i^{1+d /\bar p_1}d_j^{-d -2+d/p_1}d_j^{d/2+1}\\[5pt]
&= C(d_i/d_j)^{1+d/2-d/p_1} .
\end{align*}
If $|i-j|\leq 1$, then by applying the standard energy estimate we get
$\vertiii{w}_{2,Q_T}\leq C\vertiii{v}_{Q_T}=C$.
Combining the three cases, we have proved
$$\vertiii{w}_{2,Q_i'}\leq C\min\big(d_i/d_j,d_j/d_i\big)^{1+d/2-d/p_1}:=C m_{ij}.$$
Substituting \refe{f69}-\refe{sd80} into (\ref{dka6}) gives 
\begin{align}
\vertiii{F}_{Q_j'}\leq Ch^2d_j^{-d/2-1}
(h/d_j)^{-d/p_1}+C\sum_{*,i}m_{ij}
(h^2\vertiii{F_t}_{Q_i}
+h\vertiii{F}_{1,Q_i}) .
\end{align}
By noting that $p_1>d$, \refe{dlj6} reduces to
\begin{align*}
{\cal K}&\leq C+C\sum_j (h/d_j)^{1-d/p_1}+C\sum_j d_j^{d/2}
\mu_j \sum_{*,i}m_{ij}\big(h^2\vertiii{F_t}_{Q_i}
+h\vertiii{F}_{1,Q_i} \big)\\
&\leq
C+C\sum_{*,i}\left(h^2\vertiii{F_t}_{Q_i}+h\vertiii{F}_{1,Q_i}
\right)\sum_jd^{d/2}_j\mu_jm_{ij}\\
&\leq
C+C\sum_{*,i}\left(h^2\vertiii{F_t}_{Q_i}+h\vertiii{F}_{1,Q_i}
\right)
d^{1+d/2}_i\mu_id^{-1}_i\\
&\leq C+\left(h\vertiii{F_t}_{Q_*}
+\vertiii{F}_{1,Q_*} \right)d_{J_*}^{d/2}
\big(1/\ln(2+1/h)+h/d_{J_*}\big)\\
&~~~+C\sum_id^{1+d/2}_i
\left(\vertiii{F_t}_{Q_i}+\mu_i\vertiii{F}_{1,Q_i}\right
)\left(\frac{h}{d_i}\right)\\
&\leq
C+CC^{d/2}_*+C\sum_id^{1+d/2}_iK_i\left(\frac{h}{d_i}\right)\\
&\leq C_2+C_2C^{d/2}_*+C_2C^{-1}_*{\cal K}
\end{align*}
for some positive constant $C_2$.  
By choosing $C_*=16+2C_2$, the above inequality shows
that ${\cal K}$ is bounded. As we have mentioned,
the boundedness of ${\cal K}$ implies (\ref{FFEst1}).

Next, we prove \refe{FFEst2}. From the definition of ${\cal K}$ in (\ref{KdKj0})-(\ref{KdKj}),
we further derive that
\begin{align*}
\vertiii{\partial_tF}_{L^2(\Omega\times(1/4,1))}+
 \vertiii{\partial_{tt}F}_{L^2(\Omega\times(1/4,1))} \leq C .
\end{align*}
The above inequality and (\ref{dlk80}) imply that
\begin{align*}
\vertiii{\partial_t\Gamma_h}_{L^2(\Omega\times(1/4,1))}+
 \vertiii{\partial_{tt}\Gamma_h}_{L^2(\Omega\times(1/4,1))} \leq C ,
\end{align*}
which together with (\ref{dlky60}) gives
\begin{align*}
\|\partial_t\Gamma_h(1,\cdot,x_0)\|_{L^2}+
\|\partial_{tt}\Gamma_h(1,\cdot,x_0)\|_{L^2} \leq C .
\end{align*}
Differentiating the equation (\ref{GMhFdef}) with respect to $t$ and
multiplying the result by $\partial_t\Gamma_h$, we get
\begin{align*}
&\frac{\d}{\d t}\|\partial_t\Gamma_h(t,\cdot,x_0)\|_{L^2}^2
+c_0\|\partial_t\Gamma_h(t,\cdot,x_0)\|_{L^2}^2  \\
&\leq \frac{\d}{\d t}\|\partial_t\Gamma_h(t,\cdot,x_0)\|_{L^2}^2
+(A_h\partial_t\Gamma_h(t,\cdot,x_0),\partial_t\Gamma_h(t,\cdot,x_0)) \\
&= 0
\end{align*}
for $t\geq 1$, which further gives
\begin{align*}
\|\partial_t\Gamma_h(t,\cdot,x_0)\|_{L^2}^2\leq
e^{-c_0(t-1)}\|\partial_t\Gamma_h(1,\cdot,x_0)\|_{L^2}^2
\leq Ce^{-c_0(t-1)} .
\end{align*}
Similarly, we can prove that
\begin{align*}
\|\partial_{tt}\Gamma_h(t,\cdot,x_0)\|_{L^2}^2
\leq Ce^{-c_0(t-1)} .
\end{align*}
From (\ref{FFEst1}), (\ref{GFest0423})
and the last two inequalities, we derive (\ref{FFEst2}) for the case
$h<h_0:=1/(4C_*)$.

The proof of Lemma \ref{GMhEst} is completed. ~\endproof

\section{Proof of Theorem \ref{MainTHM1}}
\label{pkq}
\setcounter{equation}{0}
In this section, we assume that $C_*=16+2C_2$ as chosen in the last section.

\subsection{Proof of (\ref{STLEst})-(\ref{STLEst2})}
Since
$$(E_h(t)v_h)(x_0) = (F (t), v_h) + (\Gamma(t), v_h)
=\int_0^t(\partial_tF (s), v_h)\d s+(F(0), v_h) + (\Gamma(t), v_h)$$
with $\|F(0)\|_{L^1}+\|\Gamma(t)\|_{L^1}\leq C$, it follows that
(\ref{STLEst}) is a consequence of (\ref{FFEst2}) when $h<h_0$.

From Page 1360 of \cite{STW2} we also see that, for $t\in(0,T)$,
\begin{align*}
u_h(t)(x_0) = (u(t), \widetilde \delta_h) + 
\int_0^t(u(t), \partial_tF(t-s)) \d s
+\int_0^t (u(s),AF(t-s))\d s ,
\end{align*}
where the first two terms on the right-hand side
are bounded by $C\|u\|_{L^\infty(Q_T)}$ and
the third term satisfies that
\begin{align*}
&\bigg|\int_0^t (u(s),AF(t-s))\d s \bigg|\\
&\leq C\|u\|_{L^\infty(Q_T)}
\big(h^{-1}\|F\|_{W^{1,0}_1( Q_T)}+\|I_h\Gamma-\Gamma\|_{W^{2,0}_1(Q_T)}^{(h)}
+h^{-1}\|I_h\Gamma-\Gamma\|_{W^{1,0}_1(Q_T)}\big)\\
&\leq C\|u\|_{L^\infty(Q_T)}
\big(l_h+\|I_h\Gamma-\Gamma\|_{W^{2,0}_1(Q_T)}^{(h)}
+h^{-1}\|I_h\Gamma-\Gamma\|_{W^{1,0}_1(Q_T)}\big) ,
\end{align*}
where we have used (\ref{FFEst1}) in the last inequality. We
see that (\ref{STLEst2}) is a consequence of the following inequality:
\begin{align}\label{dso7}
\|I_h\Gamma-\Gamma\|_{W^{2,0}_1(Q_T)}^{(h)}
+h^{-1}\|I_h\Gamma-\Gamma\|_{W^{1,0}_1(Q_T)}
\leq Cl_h .
\end{align}
To check the above inequality, we simply note that
\begin{align*}
\|I_h\Gamma-\Gamma\|_{W^{2,0}_1(Q_*)}^{(h)}
+h^{-1}\|\nabla(I_h\Gamma-\Gamma)\|_{L^1(Q_*)}^{(h)}
&\leq
CC_*^{1+d/2}h^{1+d/2}\vertiii{\Gamma}_{2,Q_T}\\
&\leq
CC_*^{1+d/2}h^{1+d/2}\|\widetilde\delta_{x_0}\|_{H^1(\Omega)}\\
&\leq CC_*^{1+d/2},
\end{align*}
and by Lemma \ref{GFEst1} we have
\begin{align*}
&\|I_h\Gamma-\Gamma\|_{W^{2,0}_1(Q_T\backslash
Q_*)}^{(h)}+h^{-1}\|I_h\Gamma-\Gamma\|_{W^{1,0}_1(Q_T\backslash
Q_*)} \\[8pt]
&\leq C\sum_j d_j^{1+d/2}\vertiii{\Gamma}_{2,Q_j} \leq C\sum_j
d_j^{1+d/2}d_j^{-1-d/2} \\
&\leq CJ_* \leq Cl_h .
\end{align*}
Therefore, (\ref{dso7}) is proved for $T=1$, which implies (\ref{STLEst2}) for $T=1$
and $h<h_0$.
The case $T>1$ follows from the case $T=1$ by iterations:
\begin{align*}
\|u_h\|_{L^\infty(\Omega\times(k,k+1])}\leq
C\|u_h\|_{L^\infty(\Omega\times(k-1,k])}+Cl_h\|u\|_{L^\infty(\Omega\times(0,T))},
\quad\forall~k\geq 1.
\end{align*}

When $h\geq h_0$ and $f\equiv g_j\equiv 0$, the standard energy
estimates of (\ref{FEEq0}) give
\begin{align*}
\|u_h(t)\|_{L^2}+t\|\partial_tu_h(t)\|_{L^2}\leq C\|u_h(0)\|_{L^2} .
\end{align*}
By using an inverse inequality, we further derive that
\begin{align*}
\|u_h(t)\|_{L^\infty}+t\|\partial_tu_h(t)\|_{L^\infty}
\leq Ch_0^{-d/2}(\|u_h(t)\|_{L^2}+t\|\partial_tu_h(t)\|_{L^2})
\leq C\|u_h(0)\|_{L^2} \leq C\|u_h(0)\|_{L^\infty},
\end{align*}
which implies (\ref{STLEst}).

When $h\geq h_0$ while $f$ or $g_j$ may not be identically zero, we decompose
the solution of (\ref{FEEq0})  as
$u_h=\widetilde u_h+v_h$, where $\widetilde u_h$ and $v_h$ are solutions of
the equations
\begin{align}\label{Ga513}
&\left\{
\begin{array}{ll}
\partial_t\widetilde u_h+A_h\widetilde u_h = f_h-\overline\nabla_h\cdot{\bf g} ,
\\[3pt]
\widetilde u_h(0)=P_hu^0 ,
\end{array}
\right.
\end{align}
and
\begin{align}\label{Eqv9}
&\!\!\!\!\!\!\!\!\!\!\!\!\!\!\!\!\!\!\!\!\!\left\{
\begin{array}{ll}
\partial_tv_h+A_hv_h = 0 ,
\\[3pt]
v_h(0)=u^0_h-P_hu^0,
\end{array}
\right.
\end{align}
respectively.  Write the equation (\ref{PDE0}) as
\begin{align}\label{Ga903}
&\left\{
\begin{array}{ll}
\partial_tu+Au
=f -\overline\nabla\cdot {\bf g} &\mbox{in}~\Omega,
\\[3pt]
u(0)=u^0 &\mbox{in}~\Omega ,
\end{array}
\right.
\end{align}
and let $w_h=\widetilde u_h-P_hu$. The difference of (\ref{Ga513}) and
(\ref{Ga903}) gives
\begin{align*}
&\left\{
\begin{array}{ll}
\partial_t w_h+A_h w_h =   A_h(R_hu-P_hu) ,
\\[3pt]
w_h(0)=0 .
\end{array}
\right.
\end{align*}
Multiplying the above equation by $w_h$, we obtain
\begin{align*}
\|w_h\|_{L^\infty((0,T);L^2)}
&\leq C\|R_hu-P_hu\|_{L^2((0,T);H^1)}\\
&\leq Ch_0^{-1}\|R_hu-P_hu\|_{L^2((0,T);L^2)}\\
&\leq C_{T}\|R_hu-P_hu\|_{L^\infty((0,T);L^\infty)}\\
&\leq C_{T}\|u\|_{L^\infty( Q_T)} ,
\end{align*}
where we have used the inequality 
$\|R_hu\|_{L^\infty}\leq C_{h_0}\|u\|_{L^\infty}$
in the last step. By using an inverse inequality we further derive that
\begin{align*}
\|w_h\|_{L^\infty( Q_T)}\leq
Ch_0^{-d/2}\|w_h\|_{L^\infty((0,T);L^2)}\leq C_T\|u\|_{L^\infty( Q_T)}.
\end{align*}

Applying (\ref{STLEst}) to the equation (\ref{Eqv9}) we obtain
$$
\|v_h\|_{L^\infty( Q_T)}\leq C\|u_h^0-P_hu^0\|_{L^\infty}
\leq C\|u_h^0\|_{L^\infty}+C\|u\|_{L^\infty( Q_T)}  .
$$
The last two inequalities imply (\ref{STLEst2}) for the case $h\geq h_0$.

\subsection{Proof of (\ref{smgest})}

We define the truncated Green function $G_{\rm tr}^*$ in the
following way. Let $\eta$ be a nonnegative smooth function on $\R$
such that $\eta(\rho)=0$ for $|\rho|\leq 1/2$ and $\eta(\rho)=1$ for
$|\rho|\geq 1$. If we set $\chi(t,x,y)=\eta\big(|x-y|^4+t^2\big)$
and $\chi_\epsilon(t,x,y)=
\chi(t/\epsilon^2,x/\epsilon,y/\epsilon)$, then $\chi_\epsilon$ is a
$C^\infty$ function of $x,y$ and $t$. It is easy to see that $\chi_\epsilon=0$ when
$\max(|x-y|,\sqrt{t})<\epsilon/2$, and $\chi_\epsilon=1$ when
$\max(|x-y|,\sqrt{t})>\epsilon$, and $|\partial^{\alpha_1}_t
\partial^{\beta_1}_x\partial^{\beta_2}_y
\chi_\epsilon(t,x,y)|\leq C\epsilon^{-2\alpha_1
-|\beta_1|-|\beta_2|}$.

For  $d_{J_*}=C_*h$,
$\chi_{d_*}(\cdot,\cdot,y)=0$ in the domain
$Q_{*/2}(y):=\{(x,t)\in Q_T: \max(|x-y|,\sqrt{t})<d_{J_*}/2\}$.
We define a truncated Green's function by
\begin{align}
G_{\rm tr}^*(t,x,y)=
G(t,x,y)\chi_{d_{J_*}}(t,x,y) .
\end{align}
Check that $G_{\rm tr}^{*}(t,x,y)$ is symmetric with respect to $x$
and $y$, $G_{\rm tr}^{*}(\cdot,\cdot,y)\equiv 0$ in $Q_{*/2}(y)$, $0\leq G^*_{\rm
tr}(t,x,y)\leq G(t,x,y)$ and it obeys (\ref{FEstP})-(\ref{FtEstP}) when
$\max(|x-y|,\sqrt{t})>d_{J_*}$.

For the fixed trianglular element $\tau_l^h$ and
the point $x_0\in \tau_l^h$, the function
$\widetilde \delta_{x_0}$ is supported in
$\tau_l^h\subset \Omega_*(x_0)$ with
$\int_\Omega\widetilde\delta_{x_0}(y)\d y=1$
(see the notations in Section \ref{fnot}). Therefore, by using
Lemma \ref{GFEst1}  we see that
\begin{align*}
&\iint_{\Omega_\infty\backslash
Q_*(x_0)}|\partial_t\Gamma(\tau,x,x_0)
-\partial_tG(\tau,x,x_0)|\d x\d\tau\\
&= \iint_{ Q_T\backslash Q_*(x_0)}\biggl|\int_\Omega
\partial_tG(\tau,x,y)\widetilde\delta_{x_0}(y)\d
y-\partial_tG(\tau,x,x_0)\biggl|\d x\d\tau\\
&~~~
+\iint_{\Omega\times(T,\infty)}|\partial_t\Gamma(\tau,x,x_0)
-\partial_tG(\tau,x,x_0)|\d x\d\tau\\
&\leq
Ch\iint_{\max(|x-y|,\tau^{1/2})>\frac{1}{2}C_*h}\sup_{y\in\tau^l}\big|
\nabla_y\partial_tG(\tau,x,y)\big|\d x\d\tau
+C   \\
&=Ch\sum_{j}\iint_{Q_j'(y)}\sup_{y\in\tau^l}\big|
\nabla_y\partial_tG(\tau,x,y)\big|\d x\d\tau
+C\\
&\leq C\sum_{j}\frac{h}{d_j} +C \\
&\leq C .
\end{align*}

Multiplying (\ref{GMFdef}) by $\partial_t\Gamma$ and integrating the
result, we get
\begin{align*}
&\|\partial_t\Gamma(\cdot,\cdot,x_0)\|_{L^2(Q_T)}\leq
C\|\widetilde \delta_{x_0}\|_{H^1(\Omega)}\leq Ch^{-d/2-1} ,
\end{align*}
which implies that
\begin{align*}
&\iint_{Q_*(x_0)}|\partial_t\Gamma(\tau,x,x_0)|\d x\d\tau
\leq
d_{J_*}^{d/2+1}\|\partial_t\Gamma(\cdot,\cdot,x_0)\|_{L^2(Q_T)}\leq
C .
\end{align*}
Easy to check that
\begin{align}
|\partial_tG_{\rm tr}^*(t,x,y)|\leq Cd_{J_*}^{-d-2}
\quad\mbox{for}~\max(|x-y|,t^{1/2}) \leq d_{J_*}
\end{align}
and  so
\begin{align}
&\iint_{Q_*(x_0)} |\partial_tG_{\rm tr}^*(t,x,x_0)|\d
x\d\tau\leq Cd_{J_*}^{-d-2} d_{J_*}^{d+2} \leq C .
\end{align}
It follows that
\begin{align*}
&\iint_{\Omega_\infty}|\partial_t\Gamma(\tau,x,x_0)
-\partial_tG_{\rm tr}(\tau,x,x_0)|\d x\d\tau\\
&= \iint_{\Omega_\infty\backslash
Q_*}|\partial_t\Gamma(\tau,x,x_0)
-\partial_tG(\tau,x,x_0)|\d x\d\tau
+\iint_{Q_*}(|\partial_t\Gamma(\tau,x,x_0)|+
|\partial_tG_{\rm tr}(\tau,x,x_0)|)\d x\d\tau\\
&\leq C .
\end{align*}
From Lemma \ref{GMhEst} and the last inequality, we see that
\begin{align*}
&\iint_{\Omega_\infty}|\partial_t\Gamma_{h}(\tau,x,x_0)
-\partial_tG_{\rm tr}(\tau,x,x_0)|\d x\d\tau\\
&\leq\iint_{\Omega_\infty}|\partial_t\Gamma_{h}(\tau,x,x_0)
-\partial_t\Gamma(\tau,x,x_0)|\d x\d\tau
+\iint_{\Omega_\infty}|\partial_t\Gamma(\tau,x,x_0)
-\partial_tG_{\rm tr}(\tau,x,x_0)|\d x\d\tau\\
&\leq C.
\end{align*}

Since both $\Gamma_h(\tau,x,y)$ and $G^*_{\rm tr}(\tau,x,y)$
are symmetric with respect to $x$ and $y$, from the last
inequality we see that the kernel
$K(x,y)=\int_0^\infty|\partial_t\Gamma_h(\tau,x,y)
-\partial_tG^*_{\rm tr}(\tau,x,y)|\d\tau$
satisfies
\begin{align*}
&\sup_{y\in\Omega}\int_\Omega K(x,y)\d
x+\sup_{x\in\Omega}\int_\Omega K(x,y)\d y\leq C .
\end{align*}
By Schur's lemma \cite{Kra}, the operator $M_K$ defined by
$M_Ku_h(x)=\int_\Omega K(x,y)u_h(y)\d y$ is bounded on
$L^q(\Omega)$ for any $1\leq q\leq\infty$, i.e.
\begin{align}\label{Ms1}
\|M_Ku_h\|_{L^q}\leq C\|u_h\|_{L^q} ,\quad 1\leq q\leq \infty
.
\end{align}

Let
$E^*_{\rm tr}(t)u_h(x)=\int_\Omega G_{\rm tr}^*(t,x,y)u_h(y)\d y$.
We have
\begin{align*}
&\sup_{t>0}|E_h(t)u_h(x)|\\
&
\leq\sup_{t>0}|(E_h(t)-E^*_{\rm
tr}(t))u_h(x)|+\sup_{t>0}|E^*_{\rm tr}(t)u_h(x)|\\
&
\leq
|(P_h\delta_{x},u_h)|
+\sup_{t>0}\biggl|\int_0^t
\int_\Omega(
\partial_t\Gamma_h(\tau,\cdot,y)
-\partial_tG^*_{\rm tr}(\tau,x,y)u_h(y))\d
y\d\tau\biggl|+\sup_{t>0}|E^*_{\rm tr}(t)u_h(x)|\\
&\leq
|(P_h\delta_{x},u_h)|
+\int_0^\infty
\int_\Omega
|\partial_t\Gamma_h(\tau,x,y)
-\partial_tG^*_{\rm tr}(\tau,x,y)|u_h(y)\d y\d
\tau+\sup_{t>0}E(t)|u_h|(x),\\
&
=:|u_h(x)|+M_Ku_h(x)+\sup_{t>0}E(t)|u_h|(x)
\end{align*}
where
\begin{align*}
&\|M_Ku_h\|_{L^q}\leq C\|u_h\|_{L^q} , \qquad\qquad\forall~
1\leq q\leq \infty,~~\mbox{by (\ref{Ms1})},\\[5pt]
&\|\sup_{t>0}E(t)|u_h|\|_{L^q}\leq C_q\|u_h\|_{L^q},
\quad\forall~1<q<\infty ,~~\mbox{by \cite{Gra}} ,\\
&\|\sup_{t>0}E(t)|u_h|\|_{L^\infty}\leq \|u_h\|_{L^\infty},
\quad~ \mbox{by the maximum principle} .
\end{align*}
This proves (\ref{smgest}) for the case $h<h_0$.

On the other hand, when $h\geq h_0$,  from (\ref{STLEst}) we see that
\begin{align*}
\|\sup_{t>0}|E_h(t)v_h|\big\|_{L^q}
\leq C\sup_{t>0}\|E_h(t)v_h\|_{L^\infty}
\leq C\|v_h\|_{L^\infty}\leq Ch_0^{-d/q}\|v_h\|_{L^q} .
\end{align*}

The proof of (\ref{smgest}) is completed.

\subsection{Proof of (\ref{LpqSt1})-(\ref{LpqSt3})}

Since the operator $E_h(t)$ is symmetric, i.e.
$(E_h(t)u_h,v_h)=(u_h,E_h(t)v_h)$ for any $u_h,v_h\in S_h$,
from (\ref{STLEst}) we derive that, by a
duality argument and by interpolation \cite{BL},
\begin{align}
&\|E_h(t)v_h\|_{L^q} + t\|\partial_tE_h(t)v_h\|_{L^q} \leq
C\|v_h\|_{L^q}, \quad\mbox{for}~~ 1\leq q\leq\infty , \label{anf3}
\end{align}
which means that $\{E_h(t)\}_{t>0}$ is an analytic semigroup on $L^q_h$.

First, we prove (\ref{LpqSt3}).
For the case $u^0_h\equiv {\bf g}\equiv 0$,
we rewrite the equation (\ref{FEEq0}) as
\begin{align}\label{Ga5}
&\left\{
\begin{array}{ll}
\partial_tu_h+A_hu_h=f_h  ,
\\[5pt]
u_h(0)=0 ,
\end{array}
\right.
\end{align}
where $f_h=P_hf$.
From \cite{Weis1,Weis2}, we know that the maximal $L^p$
regularity (\ref{LpqSt3}) holds iff one of the following sets is
$R$-bounded in ${\cal L}(L^q_h,L^q_h)$ independent of $h$:\\
(i)~ $\{\lambda(\lambda+A_h)^{-1}:|{\rm arg}(\lambda)|<
\pi/2+\theta\}$ for some
$0<\theta< \pi/2 \,$,\\
(ii)~ $\{E_h(t),~tA_hE_h(t):t>0\} \,$, \\
(iii)~ $\{E_h(z):|{\rm arg}(z)|<\theta\}$ for some $0<\theta< \pi/2
\,$.

Moreover, from Lemma 4.c in \cite{Weis2} we know that the set in (iii) is $R$-bounded in
${\cal L}(L^q_h,L^q_h)$ for some $\theta=\theta_{\kappa_q}>0$ if the
analytic semigroup $\{E_h(z)\}$ satisfies the maximal estimate:
\begin{align*}
\biggl\|\sup_{t>0}\biggl|\frac{1}{t} \int_0^tE_h(s)u_h\d
s\biggl|\biggl\|_{L^q}\leq \kappa_q\|u_h\|_{L^q},\quad\forall~u_h\in
L^q_h(\Omega) .
\end{align*}
Since the last inequality is a consequence of the maximal semigroup estimate
(\ref{smgest}), we thus proved the maximal $L^p$ regularity
(\ref{LpqSt3}).

Secondly, we prove \refe{LpqSt1} and \refe{LpqSt2}. For the general case $u^0_h\neq 0$
or ${\bf g}\neq 0$, we let $u_h=\widetilde u_h+v_h$, where
$\widetilde u_h$ and $v_h$ are the solutions of the equations
\begin{align}\label{Ga51}
&\left\{
\begin{array}{ll}
\partial_t\widetilde u_h+A_h\widetilde u_h = f_h-\overline\nabla_h\cdot{\bf g} ,
\\[3pt]
\widetilde u_h(0)=P_hu^0 ,
\end{array}
\right.
\end{align}
and
\begin{align}
&\!\!\!\!\!\!\!\!\!\!\!\!\!\!\!\!\!\!\!\!\!\left\{
\begin{array}{ll}
\partial_tv_h+A_hv_h = 0 ,
\\[3pt]
v_h(0)=u^0_h-P_hu^0,
\end{array}
\right.
\end{align}
respectively.  Write the equation (\ref{PDE0}) as
\begin{align}\label{Ga90}
&\left\{
\begin{array}{ll}
\partial_tu+Au
=f -\overline\nabla\cdot {\bf g} &\mbox{in}~\Omega,
\\[3pt]
u(0)=u^0 &\mbox{in}~\Omega ,
\end{array}
\right.
\end{align}
and let $w_h=\widetilde u_h-P_hu$. The difference of (\ref{Ga51}) and
(\ref{Ga90}) gives
\begin{align*}
&\left\{
\begin{array}{ll}
\partial_t w_h+A_h w_h =   A_h(R_hu-P_hu) ,
\\[3pt]
w_h(0)=0 .
\end{array}
\right.
\end{align*}
Multiplying the above equation by $A_h^{-1}$, we get
\begin{align*}
&\left\{
\begin{array}{ll}
\partial_t A_h^{-1}w_h+A_h A_h^{-1}w_h = R_hu-P_hu ,
\\[3pt]
A_h^{-1}w_h(0)=0 ,
\end{array}
\right.
\end{align*}
and using (\ref{LpqSt3}) we derive that
\begin{align*}
\|w_h\|_{L^p((0,T);L^q)}\leq C_{p,q}\|R_hu-P_hu \|_{L^p((0,T);L^q)} .
\end{align*}
On the other hand, it is easy to derive that $\|v_h(t)\|_{L^2}\leq
Ce^{-t/C}\|u^0_h-P_hu^0\|_{L^2}$, which with (\ref{anf3})
gives (via interpolation)
\begin{align*} \|v_h(t)\|_{L^q}\leq
Ce^{-t/C_q}\|u^0_h-P_hu^0\|_{L^q} ~~~\mbox{for}~~1< q<\infty .
\end{align*}
The last two inequalities imply (\ref{LpqSt1}).

If $u^0_h\equiv u^0\equiv f\equiv 0$, then $v_h=0$
and by using Lemma \ref{s00} we derive that
\begin{align*}
\|u_h\|_{L^p((0,T);W^{1,q})}
&=\|\widetilde u_h\|_{L^p((0,T);W^{1,q})} \\
&\leq \|w_h\|_{L^p((0,T);W^{1,q})}
+\|P_hu\|_{L^p((0,T);W^{1,q})}\\
&\leq Ch^{-1}\|w_h\|_{L^p((0,T);L^q)}
+C\|u\|_{L^p((0,T);W^{1,q})}\\
&\leq C_{p,q}h^{-1}\|R_hu-P_hu \|_{L^p((0,T);L^q)}
+C\|u\|_{L^p((0,T);W^{1,q})} \\
&\leq C_{p,q}\|u\|_{L^p((0,T);W^{1,q})} \\
&\leq C_{p,q}\|{\bf g}\|_{L^p((0,T);L^q)} .
\end{align*}
This proves the inequality (\ref{LpqSt2}).

The proof of Theorem \ref{MainTHM1} is completed. ~\endproof

\section*{Appendix --- 
Interior estimates of parabolic equations on $\bf\widetilde Q_j'$}
\renewcommand{\thelemma}{A.\arabic{lemma}}
\renewcommand{\theproposition}{A.\arabic{lemma}}
\renewcommand{\theequation}{A.\arabic{equation}}
\setcounter{lemma}{0} \setcounter{equation}{0}

\appendix
\addcontentsline{toc}{section}{Appendix}
\addtocontents{toc}{\protect\setcounter{tocdepth}{2}}

\begin{lemma}\label{LemAP1}
{\it Suppose that $\widetilde\phi$ is the solution of 
\begin{align}\label{Eqofphi}
\partial_t\widetilde\phi-\nabla_{\widetilde x}\cdot(
\widetilde a\nabla_{\widetilde x}\widetilde\phi)
+\widetilde c\widetilde\phi=0\quad \mbox{in}~~\widetilde Q_j',
\end{align}
with the Neumann boundary condition 
$\widetilde a\nabla_{\widetilde x}\widetilde\phi \cdot{\bf n}=0$
on $\widetilde\Omega_j'\cap\partial\widetilde\Omega$,
then we have
\begin{align*}
&|\!|\!|\partial_{\widetilde t}\widetilde
\phi|\!|\!|_{\widetilde Q_j}+|\!|\!|\widetilde
\phi|\!|\!|_{2,\widetilde Q_j}
+|\!|\!|\partial_{\widetilde t}\widetilde
\phi|\!|\!|_{2,\widetilde Q_j}
+|\!|\!|\partial_{\widetilde t\widetilde t}\widetilde
\phi|\!|\!|_{2,\widetilde Q_j}
\leq  C|\!|\!|\widetilde\phi|\!|\!|_{\widetilde Q_j'} ,\\
&\|\partial_{\widetilde x_i}\widetilde
\phi\|_{L^\infty(\widetilde Q_j)}
+\|\partial_{\widetilde x_i}\widetilde
\phi\|_{C^{\alpha,\alpha/2}(\overline{\widetilde Q}_j)}
+ \|\partial_{\widetilde x_i}\partial_{\widetilde x_l}\widetilde
\phi\|_{L^{\infty,p_1}(\widetilde Q_j)} 
\leq C|\!|\!|\widetilde\phi|\!|\!|_{\widetilde Q_j'} . 
\end{align*}
}
\end{lemma}
\noindent {\it Proof}~~~
Without loss of generality, we can assume that
the equation (\ref{Eqofphi}) holds in $\widetilde Q_j'''$
with the boundary condition 
on $\widetilde\Omega_j'''\cap\partial\widetilde\Omega$,
and then kick $\widetilde Q_j'''$ back to $\widetilde Q_j'$
to draw the conclusion.

Let $\widetilde\omega(x,t)$ be a smooth function which has compact 
support in
$\widetilde Q_j'$ and equal to $1$ in $\widetilde Q_{j}$,
with $|\nabla \widetilde\omega|\leq C$ and 
$|\partial_t\widetilde\omega|\le C$,
and let $\widetilde\chi(x,t)$ be a smooth function which has compact 
support in
$\widetilde Q_j''$ and equal to $1$ in $\widetilde Q_{j}'$,
with $|\nabla \widetilde\chi|\leq C$ and 
$|\partial_t\widetilde\chi|\le C$.
Since $ \cup_{k\leq j}\widetilde\Omega_k''
\cup \widetilde\Omega_*$ is of unit size,
there exists a smooth subdomain 
$\widetilde D\subset \widetilde\Omega$
such that $\widetilde D$ has unit size and contains 
$ \cup_{k\leq j}\widetilde\Omega_k''
\cup \widetilde\Omega_*$ ($\partial\widetilde D$ may contain
a part of the boundary of $\partial\widetilde\Omega$). Then
$\widetilde D\times(0,16)$ contains $\widetilde Q_j''$ and
$\widetilde\omega \widetilde G$ 
vanishes outside $\widetilde D\times(0,16)$.
Thus $\widetilde\chi=1$ on the support
of $\widetilde\omega$.

Integrating the equation (\ref{Eqofphi}) against 
$\widetilde\omega^2\widetilde\phi$, we derive the basic local energy estimate
\begin{align}\label{nbyGq}
\|\partial_{\widetilde t}(\widetilde\omega \widetilde\phi)\|_{L^2((0,16);H^{-1}(\widetilde D))}
+\|\nabla_{\widetilde x}(\widetilde\omega \widetilde\phi)\|_{L^{2}(\widetilde D\times(0,16))} \leq 
C|\!|\!|\widetilde\phi|\!|\!|_{\widetilde Q_j'} .
\end{align}
To present further estimates for $\widetilde\phi$,
we consider $\widetilde \omega\widetilde\phi$, which is 
the solution of 
\begin{align*}
&\left\{\begin{array}{ll}
\partial_{\widetilde t}(\widetilde\omega\widetilde\phi)
-\nabla_{\widetilde x} \cdot(\widetilde a 
\nabla_{\widetilde x} (\widetilde\omega\widetilde\phi))
=\widetilde\chi\widetilde\phi\partial_{\widetilde t}\widetilde\omega-\widetilde
a\nabla_{\widetilde x} \widetilde\omega\cdot\nabla_{\widetilde x}
(\widetilde\chi\widetilde\phi)
-\nabla_{\widetilde x} \cdot(\widetilde 
a\widetilde\chi\widetilde\phi\nabla_{\widetilde x} \widetilde\omega)
& \mbox{in~ $\widetilde D\times(0,16)$} ,\\[8pt]
\widetilde a 
\nabla_{\widetilde x} (\widetilde\omega\widetilde\phi)\cdot\widetilde{\bf n}=
\widetilde 
a\widetilde\chi\widetilde\phi\nabla_{\widetilde x} \widetilde\omega\cdot\widetilde{\bf n}
&\mbox{on~ $\partial\widetilde D\times(0,16)$} ,\\[8pt]
\widetilde\omega\widetilde\phi=0
&\mbox{on~ $\widetilde D\times\{0\}$} .
\end{array}\right.
\end{align*}
Since $\widetilde D$ is a regular domain, 
the classical $L^p((0,16);W^{1,p}(\widetilde\Omega))$ 
estimate (see Lemma 2.1 and \refe{dsk4})  gives 
\begin{align*}
&\|\partial_{\widetilde t}(\omega\widetilde\phi)\|_{L^{p}((0,16);
W^{-1,p}(\widetilde D))}+\|\nabla_{\widetilde x}(\widetilde\omega
\widetilde\phi)\|_{L^{p}(\widetilde D\times(0,16))}
+\|\widetilde\omega\widetilde\phi\|_{L^{p}(\widetilde D\times(0,16))}\\
&\leq 
C\|\widetilde \chi\widetilde\phi  \|_{L^{p}(\widetilde D\times(0,16))}
+C\| \nabla_{\widetilde x}(\widetilde \chi\widetilde\phi) 
 \|_{L^{q}(\widetilde D\times(0,16))} \\
&\leq C\|\partial_{\widetilde t}(\widetilde\chi\widetilde\phi)
\|_{L^{q}((0,16);
W^{-1,q}(\widetilde D))}+
C\|\nabla_{\widetilde x}(\widetilde\chi
\widetilde\phi)\|_{L^{q}(\widetilde D\times(0,16))}
+C\|\widetilde\chi\widetilde\phi\|_{L^{q}(\widetilde D\times(0,16))}
\end{align*}
for any $p>2$ and $q=(d+2)p/(d+2+p)<p$, where
we have used the Sobolev embedding 
$L^q((0,16);W^{1,q}(\widetilde D))\cap 
W^{1,q}((0,16);W^{-1,q}(\widetilde D))\hookrightarrow
L^{p}(\widetilde D\times(0,16)) . $
Via a kickback argument and using \refe{nbyGq}, 
a finite number of iterations of the last inequality
give 
\begin{align}
\|\nabla_{\widetilde x} \widetilde\phi\|_{L^{p}(\widetilde Q_j)}
\leq  C|\!|\!|\widetilde\phi|\!|\!|_{\widetilde Q_j'} .
\label{ogpqj2}
\end{align}

Next, we present $L^p((0,16);W^{2,p}(\widetilde D))$ 
estimates for $\widetilde\phi$.
For this purpose, we let $\widetilde W$ be the solution of
\begin{align*}
&\left\{\begin{array}{ll}
\nabla_{\widetilde x} \cdot(\widetilde a 
\nabla_{\widetilde x} \widetilde W)
=\nabla_{\widetilde x} \cdot(\widetilde 
a\widetilde\chi\widetilde\phi\nabla_{\widetilde x} \widetilde\omega)
& \mbox{in~ $\widetilde D\times(0,16)$} ,\\[8pt]
\widetilde a 
\nabla_{\widetilde x} \widetilde W\cdot\widetilde{\bf n}=
\widetilde 
a\widetilde\chi\widetilde\phi\nabla_{\widetilde x} \widetilde\omega\cdot\widetilde{\bf n}
&\mbox{on~ $\partial\widetilde D\times(0,16)$} ,
\end{array}\right.
\end{align*}
and so 
$\partial_{\widetilde t}\widetilde W$ is the solution of
\begin{align*}
&\left\{\begin{array}{ll}
\nabla_{\widetilde x} \cdot(\widetilde a 
\nabla_{\widetilde x} \partial_{\widetilde t}\widetilde W)
=\nabla_{\widetilde x} \cdot(\widetilde 
a\partial_{\widetilde t}(\widetilde\chi\widetilde\phi)\nabla_{\widetilde x} \widetilde\omega
+\widetilde 
a\widetilde\chi\widetilde\phi\nabla_{\widetilde x} \partial_{\widetilde t}\widetilde\omega)
& \mbox{in~ $\widetilde D\times(0,16)$} ,\\[8pt]
\widetilde a 
\nabla_{\widetilde x} \partial_{\widetilde t}\widetilde W\cdot\widetilde{\bf n}=
(\widetilde 
a\partial_{\widetilde t}(\widetilde\chi\widetilde\phi)\nabla_{\widetilde x} \widetilde\omega
+\widetilde 
a\widetilde\chi\widetilde\phi\nabla_{\widetilde x} \partial_{\widetilde t}\widetilde\omega)
\cdot\widetilde{\bf n}
&\mbox{on~ $\partial\widetilde D\times(0,16)$} .
\end{array}\right.
\end{align*}
To estimate $\| \partial_{\widetilde t}\widetilde W
\|_{L^{p}(\widetilde D)}
$, we define $\phi$ as the solution of
\begin{align*}
&\left\{\begin{array}{ll}
-\nabla_{\widetilde x} \cdot(\widetilde a 
\nabla_{\widetilde x} \phi)
=\varphi
& \mbox{in~ $\widetilde D$} ,\\[8pt]
\widetilde a 
\nabla_{\widetilde x} \phi
\cdot\widetilde{\bf n}=0
&\mbox{on~ $\partial\widetilde D$} ,
\end{array}\right.
\end{align*}
and consider
\begin{align*}
\big|\big( \partial_{\widetilde t}\widetilde W,
\varphi\big)\big|
&= \big|\big(\widetilde a 
\nabla_{\widetilde x} \partial_{\widetilde t}\widetilde W,
\nabla_{\widetilde x} \phi\big)\big|\\
&=\big|\big(\widetilde 
a\partial_{\widetilde t}\widetilde\phi
\nabla_{\widetilde x} \widetilde\omega
+\widetilde 
a\widetilde\chi\widetilde\phi\nabla_{\widetilde x} 
\partial_{\widetilde t}\widetilde\omega,
\nabla_{\widetilde x}\phi)\big|\\
&\leq \big|\big(-\nabla_{\widetilde
x}\cdot(\widetilde a\nabla_{\widetilde x}
\widetilde\phi)+\widetilde c\widetilde\phi, 
\widetilde a\nabla_{\widetilde x}\widetilde w
\cdot\nabla_{\widetilde x}\phi\big)\big|
+ \big|\big(\widetilde 
a\widetilde\chi\widetilde\phi\nabla_{\widetilde x} 
\partial_{\widetilde t}\widetilde\omega,
\nabla_{\widetilde x}\phi)\big| \\
&\leq C\|\widetilde\phi\widetilde\chi\|_{W^{1,p}(\widetilde D)} 
\|\phi\|_{W^{2,p'}(\widetilde D)}\\
&\leq C\|\widetilde\phi\widetilde\chi\|_{W^{1,p}(\widetilde D)} 
\|\varphi\|_{L^{p'}(\widetilde D)},
\end{align*}
which implies that, via a duality argument,
\begin{align*}
\| \partial_{\widetilde t}\widetilde W
\|_{L^{p}(\widetilde D)}
\leq C\|\widetilde\chi\widetilde\phi\|_{W^{1,p}(\widetilde D)}  .
\end{align*}
With the last inequality, we note that 
$\widetilde Z=\widetilde\omega\widetilde\phi-\widetilde W$ is the
solution of
\begin{align*}
&\left\{\begin{array}{ll}
\partial_{\widetilde t}\widetilde Z
-\nabla_{\widetilde x} \cdot(\widetilde a 
\nabla_{\widetilde x} \widetilde Z)
=-\partial_{\widetilde t}\widetilde W
+\widetilde\phi\partial_{\widetilde t}\widetilde\omega-\widetilde
a\nabla_{\widetilde x} \widetilde\omega\cdot\nabla_{\widetilde x}
\widetilde\phi
& \mbox{in~ $\widetilde D\times(0,16)$} ,\\[8pt]
\widetilde a 
\nabla_{\widetilde x} \widetilde Z\cdot\widetilde{\bf n}= 0
&\mbox{on~ $\partial\widetilde D\times(0,16)$} ,\\[8pt]
\widetilde Z=0
&\mbox{on~ $\widetilde D\times\{0\}$} ,
\end{array}\right.
\end{align*}
which obeys the $L^p((0,16);W^{2,p}(\widetilde D))$ estimate  (see Lemma 2.1) 
\begin{align*}
\|\partial_{\widetilde t}\widetilde Z\|_{L^p(\widetilde D\times(0,16))}
+\|\partial_{\widetilde x_i}\partial_{\widetilde x_l}\widetilde Z\|_{L^p(\widetilde D\times(0,16))}
&\leq C\|\partial_{\widetilde t}\widetilde W\|_{L^p((0,16);L^{p}(\widetilde D))} 
+C\|\widetilde\phi\widetilde\chi\|_{L^p((0,16);W^{1,p}(\widetilde D))} \\
&\leq C\|\widetilde\phi\widetilde\chi\|_{L^p((0,16);W^{1,p}(\widetilde D))} ,
\end{align*}
and this further leads to
\begin{align}
&\|\partial_{\widetilde t}(\widetilde\omega\widetilde\phi)\|_{L^{p}(\widetilde D\times(0,16))}
+\|\partial_{\widetilde x_i}\partial_{\widetilde x_l}(\widetilde\omega\widetilde\phi)\|_{L^{p}(\widetilde D\times(0,16))}
\leq C\|\widetilde\phi\|_{L^{p}(\widetilde  Q_j'')}
+C\|\nabla_{\widetilde x}\widetilde\phi\|_{L^{p}(\widetilde  Q_j'')} .
\label{ogpqj1}
\end{align}
Then a kickback argument in 
\refe{ogpqj2}-\refe{ogpqj1} gives
\begin{align}
&\|\partial_{\widetilde t}(\widetilde\omega\widetilde\phi)\|_{L^{p}(\widetilde D\times(0,16))}
+\|\partial_{\widetilde x_i}\partial_{\widetilde x_l}(\widetilde\omega\widetilde\phi)\|_{L^{p}(\widetilde D\times(0,16))}
\leq C|\!|\!|\widetilde\phi|\!|\!|_{\widetilde Q_j'}   .
\label{ogpqjW21p}
\end{align}
Since $\partial_{\widetilde t}\widetilde\phi$, $\partial_{\widetilde t}\widetilde\phi$ and $\widetilde\phi_{y_l}$ satisfies the same equation as $\widetilde\phi$ in $Q_j'$, 
the last inequality holds if
$\widetilde\phi$ is replaced by $\partial_{\widetilde t}\widetilde\phi$, $\partial_{\widetilde t}\widetilde\phi$  or $\widetilde\phi_{y_l}$, i.e.
\begin{align}
&\|\partial_{\widetilde t}(\widetilde\omega\partial_{\widetilde t}\widetilde\phi)\|_{L^{p}(\widetilde  D\times(0,16))}
+\|\partial_{\widetilde x_i}\partial_{\widetilde x_l}(\widetilde\omega\partial_{\widetilde t}\widetilde\phi )\|_{L^{p}(\widetilde D\times(0,16))}
\leq C|\!|\!|\partial_{\widetilde t}\widetilde\phi|\!|\!|_{\widetilde Q_j'}  ,
\label{ogpqjW21p20}\\
&\|\partial_{\widetilde t}(\widetilde\omega\partial_{\widetilde t\widetilde t}\widetilde\phi)\|_{L^{p}(\widetilde  D\times(0,16))}
+\|\partial_{\widetilde x_i}\partial_{\widetilde x_l}(\widetilde\omega\partial_{\widetilde t\widetilde t}\widetilde\phi )\|_{L^{p}(\widetilde D\times(0,16))}
\leq C|\!|\!|\partial_{\widetilde t\widetilde t}\widetilde\phi|\!|\!|_{\widetilde Q_j'}  .
\label{ogpqjW21p21}
\end{align}

Via a kickback argument, 
the last three inequalities imply that
\begin{align}
&|\!|\!|\partial_{\widetilde t}\widetilde
\phi|\!|\!|_{\widetilde Q_j}+|\!|\!|\widetilde
\phi|\!|\!|_{2,\widetilde Q_j}
+|\!|\!|\partial_{\widetilde t}\widetilde
\phi|\!|\!|_{2,\widetilde Q_j}
+|\!|\!|\partial_{\widetilde t\widetilde t}\widetilde
\phi|\!|\!|_{2,\widetilde Q_j}
\leq  C|\!|\!|\widetilde\phi|\!|\!|_{\widetilde Q_j'}  ,  \\
&\|\partial_{\widetilde
x_i}\partial_{\widetilde x_l}\widetilde\phi\|_{L^{p_1}(\widetilde
Q_j)}+
\|\partial_{\widetilde
x_i}\partial_{\widetilde x_l}\partial_{\widetilde t}\widetilde\phi\|_{L^{p_1}(\widetilde
Q_j)}
+\|\widetilde\omega\widetilde\phi\|_{L^{\infty}((0,16);L^{p_1}(\widetilde D))}\leq  C|\!|\!|\widetilde\phi|\!|\!|_{\widetilde Q_j'}  , \nn
\end{align}
and so
\begin{align*}
\|\partial_{\widetilde t}\partial_{\widetilde x_i}\partial_{\widetilde x_l}(\widetilde\omega\widetilde\phi)\|_{L^{p}(\widetilde D\times(0,16))}
&\leq 
\|\partial_{\widetilde x_i}\partial_{\widetilde x_l}(\widetilde\phi\partial_{\widetilde t}\widetilde\omega)\|_{L^{p}(\widetilde D\times(0,16))}
+\|\partial_{\widetilde x_i}\partial_{\widetilde x_l}(\widetilde\omega\partial_{\widetilde t}\widetilde\phi)\|_{L^{p}(\widetilde D\times(0,16))} \\
&
\leq C\|\partial_{\widetilde
x_i}\partial_{\widetilde x_l}\widetilde\phi\|_{L^{p_1}(\widetilde
Q_j')} 
+C\|\partial_{\widetilde
x_i}\partial_{\widetilde x_l}\partial_{\widetilde t}\widetilde\phi
\|_{L^{p_1}(\widetilde
Q_j')} \\
&\leq C|\!|\!|\widetilde\phi|\!|\!|_{\widetilde Q_j''}
+C|\!|\!|\partial_{\widetilde t}\widetilde\phi|\!|\!|_{\widetilde Q_j''}  \\
&\leq C|\!|\!|\widetilde\phi|\!|\!|_{\widetilde Q_j'''} .
\end{align*}
This gives an estimate of $\widetilde\omega\widetilde\phi$ in terms of the norm of $W^{1,p_1}((0,16);W^{2,p_1}(\widetilde D))$. Since \linebreak $W^{1,p_1}((0,16);W^{2,p_1}(\widetilde D))\hookrightarrow C^{1+\alpha,(1+\alpha)/2}(\overline{\widetilde D})$ for some $\alpha>0$,  we further derive that
\begin{align*}
&\|\partial_{\widetilde x_i}(\widetilde\omega\widetilde\phi)\|_{L^{\infty}(\widetilde D\times(0,16))}
+
\|\partial_{\widetilde x_i}(\widetilde\omega\widetilde\phi)\|_{C^{\alpha,\alpha/2}(\overline{\widetilde D}\times[0,16])}
+\|\partial_{\widetilde x_i}\partial_{\widetilde x_l}(\widetilde\omega\widetilde\phi)\|_{L^{\infty}((0,16);L^{p_1}(\widetilde D))}\\
&\leq C\|\widetilde\omega\widetilde\phi\|_{L^2(\widetilde D\times(0,16))}+
C\|\partial_{\widetilde x_i}\partial_{\widetilde x_l}(\widetilde\omega\widetilde\phi)\|_{L^{p}(\widetilde D\times(0,16))}+C\|\partial_{\widetilde t}\partial_{\widetilde x_i}\partial_{\widetilde x_l}(\widetilde\omega\widetilde\phi)\|_{L^{p}(\widetilde D\times(0,16))}\\
&  \leq C|\!|\!|\widetilde\phi|\!|\!|_{\widetilde Q_j'''} ,
\end{align*}
and  this inequality also holds if $\widetilde\phi$ is replaced
by $\widetilde\phi_{y_l}$. Kicking $\widetilde Q_j'''$ back to $\widetilde Q_j'$, we derive that
\begin{align}
&\|\partial_{\widetilde x_i}\widetilde
\phi\|_{L^\infty(\widetilde Q_j)}
+\|\partial_{\widetilde x_i}\widetilde
\phi\|_{C^{\alpha,\alpha/2}(\overline{\widetilde Q}_j)}
+
\|\partial_{\widetilde x_i}\partial_{\widetilde x_l}\widetilde
\phi\|_{L^{\infty,p_1}(\widetilde Q_j)} 
\leq C|\!|\!|\widetilde\phi|\!|\!|_{\widetilde Q_j'} .
\end{align}

The proof of Lemma \ref{LemAP1} is completed. ~\endproof

\bigskip
\medskip

{\bf Acknowledgement}~~ The author would like to thank 
Professor Weiwei Sun 
for the helpful discussions
and the anonymous referees 
for the valuable comments and suggestions.

\end{document}